\documentclass[nonanon,onecolumn,11pt]{emiclass}

\usepackage{tikz}
\usepackage{subcaption}
\usepackage{xcolor}
\usepackage{graphicx}
\usepackage[colorlinks=true,linkcolor=blue,citecolor=red,urlcolor=blue,pdfborder={0 0 0}]{hyperref}

\definecolor{darkgreen}{RGB}{0,180,40}
\definecolor{darkred}{RGB}{180,0,40}

\newcommand{\st}{\textsc{st}}
\newcommand{\bo}{\textsc{bo}} 
\newcommand{\br}{\textsc{br}}
\newcommand{\brl}{\textsc{br-l}}
\newcommand{\brr}{\textsc{br-r}}
\newcommand{\su}{\textsc{su}}
\newcommand{\att}{\textsc{exc}}
\newcommand{\pl}{\textsc{pl}}
\newcommand{\bow}{\textsc{bow}}
\newcommand{\rel}{\textsc{rel}}
\newcommand{\rec}{\textsc{rec}}
\newcommand{\wolf}{\textsc{wolf}}
\newcommand{\fid}{\textsc{fidelity}}
\newcommand{\clos}{\textsc{sustain}}

\newcommand{\maxsc}{\textsc{max}}
\newcommand{\dyn}{\textsc{dyn}}
\newcommand{\stc}{\textsc{stc}}
\newcommand{\band}{\textsc{band}}
\newcommand{\tot}{\textsc{tot}}

\newcommand{\Nfreq}{N_f}
\newcommand{\Fstrbr}{F_{\st \rightleftharpoons \br}}
\newcommand{\Fbobrleft}{F_{\bo \rightleftharpoons \brl}}
\newcommand{\Fbobrright}{F_{\bo \rightleftharpoons \brr}}
\newcommand{\Fbosupp}{F_{\bo \rightleftharpoons \su}}
\newcommand{\Fattack}{F_\att}

\newcommand{\FNstrbr}{F^n_{\st \rightleftharpoons \br}}
\newcommand{\FNbobrleft}{F^n_{\bo \rightleftharpoons \brl}}
\newcommand{\FNbobrright}{F^n_{\bo \rightleftharpoons \brr}}
\newcommand{\FNbosupp}{F^n_{\bo \rightleftharpoons \su}}
\newcommand{\FNattack}{F^n_\att}

\newcommand{\audionotelocal}[2][1.35ex]{\href{#2}{
 \begin{tikzpicture}[baseline=-0.6ex]
	\draw[fill=white,  draw=black!60, line width=0.3pt] (0,0) circle (\dimexpr 3#1/25\relax);
	\node at (0,0){\includegraphics[height=\dimexpr 9#1/55\relax]{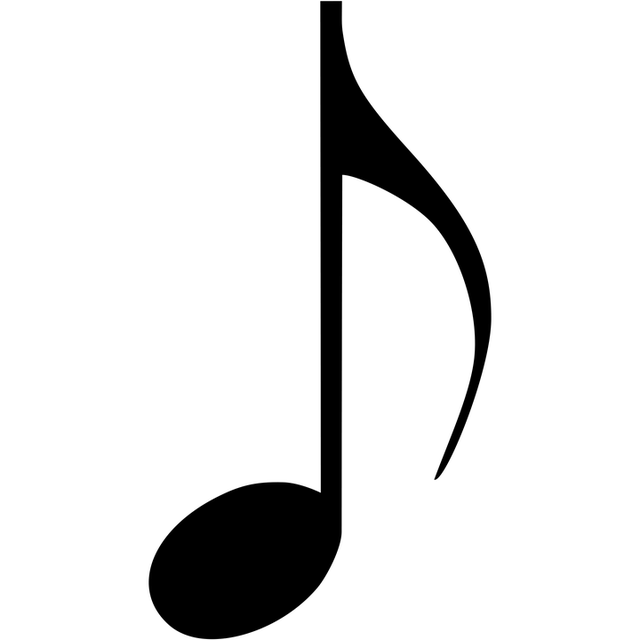}};
 \end{tikzpicture}
 }
}
\newcommand{\audionoteonline}[2][1.35ex]{\href{#2}{
 \begin{tikzpicture}[baseline=-0.6ex]
	\draw[fill=yellow, draw=black!60, line width=0.3pt] (0,0) circle (\dimexpr 3#1/25\relax);
	\node at (0,0){\includegraphics[height=\dimexpr 9#1/55\relax]{figures/croma.png}};
 \end{tikzpicture}
}
}

\title{The Wolf and the Cello: Modelling and Design of \\ Multiple Resonance Suppressors in Large String Instruments}

\author[1]{Simone Cacace}
\author[2]{Emiliano Cristiani}
\author[2,3,*]{Francesca L. Ignoto}
\affil[1]{Dipartimento di Matematica, Sapienza Università di Roma, Rome, Italy}
\affil[2]{Istituto per le Applicazioni del Calcolo, Consiglio Nazionale delle Ricerche, Rome, Italy}
\affil[3]{Dipartimento di Scienze di Base e Applicate per l'Ingegneria, Sapienza Università di Roma, Rome, Italy}
\corremail{francescalourdes.ignoto@uniroma1.it}
\date{\today}

\begin{document}
\maketitle
	
\begin{abstract}
	The wolf note is an acoustic instability that occurs in large bowed string instruments when a strong body resonance interacts with the vibrating string, producing amplitude modulation and loss of tonal control. 
	Various wolf suppressors -- tuned mass dampers attached to the string or to the instrument body -- are used in practice to mitigate this effect.
	In this paper, we propose a mathematical model describing the coupled dynamics of a string and a two-dimensional body equipped with one or two wolf suppressors. 
	Both string and body include elastic (second-order) and stiffness (fourth-order) contributions and can be excited either by plucking or bowing.
	Three performance indicators are introduced: The first one perceives the wolf-tone appearance, the second one quantifies the attenuation of the notes possibly caused by the wolf suppressor, and the third one measures the acoustic fidelity (in terms of spectrum) with respect to the original instrument. 
	The proposed numerical tests give insights about optimal tuning and placement of one or two suppressors, achieving effective wolf-note suppression while preserving as much as possible the global tonal balance.
\end{abstract}
	


\section{Introduction}\label{sec:intro} 
Large bowed string instruments exhibit a rich variety of nonlinear vibrational phenomena arising from the strong coupling between the string, the bridge, and the instrument body. While this coupling is essential for sound radiation and tonal richness, it can also give rise to undesirable instabilities that compromise playability. 
Among these, the so-called \emph{wolf note} is one of the most prominent and long-standing issues affecting instruments such as the cello and the double bass.
The wolf note manifests as a pronounced amplitude modulation which produces an uneven or ``howling'' sound, typically occurring when the frequency of a played note lies close to the frequency of a particular low-damped body mode. 
In this regime, energy is periodically exchanged between the vibrating string and the resonant body mode, leading to a loss of tonal stability and control. From the performer's perspective, this results in a note that is difficult to sustain and highly sensitive to bow pressure and speed. Also, the wolf note typically compromises the overall tonal quality. 

Describing the onset of the wolf note and developing effective strategies for its suppression is therefore a problem of both practical and scientific relevance.

\paragraph{Available resonance suppressors.}
Since the wolf note represents one of the most undesirable and frustrating defects of bowed string instruments, over the years a variety of \textit{resonance suppressors} have been developed to mitigate the problem, particularly in cellos. 
These devices, also known as \textit{wolf eliminators or terminators or killers}, are small mechanical attachments that modify the effective vibrational response of the instrument by introducing an additional local resonance. 

The underlying idea is to partially cancel the undesired coupling between the string and the body, thus stabilizing the oscillation around the problematic frequency.
Conventional wolf eliminators are based on mass dampers, which must be tuned to a fixed frequency. 
Therefore, they are only capable of eliminating that single wolf note at the specific frequency to which they are tuned \cite{neubauer2018}.
This means that, depending on their design, they may primarily dissipate energy near a target resonance or modify the coupled modal structure, shifting the system away from resonance coincidence.
Commercial suppressors typically consist of a small mass attached to a compliant element (e.g., a rubber or metal spring) mounted on the afterlength of the string (between the bridge and the tailpiece) or directly on the body. 
In the second case, they can be mounted externally on top plate, or inside the body. 
If inside, they can either be mounted permanently in a fixed position or they can be moved by an external magnet.

The effectiveness of these devices has been empirically demonstrated by luthiers and acousticians for decades, but their influence on the global tonal quality of the instrument remains a delicate trade-off: while they can attenuate the wolf note, they often alter the response at nearby frequencies, sometimes resulting in a perceptible loss of brilliance or sensitivity.
In addition, other difficulties come into play: 
First, environmental conditions (typically temperature and humidity) affect the structural properties of the instrument, such as the stiffness, thus altering the frequency of the wolf note \cite{neubauer2018}. 
Second, eliminators are difficult to place because of the high sensitivity of their effect with respect to their position.

In conclusion, the design of efficient and minimally invasive suppressors remains largely based on empirical tuning and trial-and-error procedures.

\paragraph{Relevant literature.}
Scientific literature recognized the wolf note as an issue since the early twentieth century; see, e.g., \cite{raman1916}. 
Some pioneering works like \cite{firth1973, mcintyre1978, mcintyre1979, mcintyre1983} started to reproduce this peculiar phenomenon by means of numerical simulations.
In recent years, it has been the subject of continuous experimental and theoretical research.
A representative contribution is the work \cite{inacio2008}, where the wolf note is investigated through a string-body interaction model based on modal decomposition. 
In this approach, both the string and the body are described in terms of truncated modal expansions, and the coupled dynamics is obtained by integrating a finite-dimensional system of ordinary differential equations, with body parameters identified from experimental impulse-response or admittance measurements at the bridge.

Another numerical attempt to reproduce the wolf note using highly simplified string-body coupled models can be found in \cite{ogura2010}, where beat-like amplitude modulations are obtained both in simulations and experiments.

A closely related experimental and numerical investigation is presented in the work \cite{debut2010}, which studies the effectiveness of a traditional cello wolf eliminator mounted on the string (between the bridge and the tailpiece, as usual). 
Through systematic hand-bow experiments on a real cello, supported by numerical simulations, the authors analyse how the occurrence and severity of the wolf note depend on the eliminator mass and position.
Their results show that the action of the device is highly sensitive to both parameters: Depending on the configuration, the wolf note may be attenuated, shifted in frequency, or even exacerbated. The study adopts a fully coupled modal formulation for the string, the body, and the eliminator, with model parameters identified from experimental measurements.


Another particularly relevant contribution is \cite{neubauer2018}, which combines both modelling and experimental validation of the wolf note and its suppression. 
Despite relying on an extremely simple linear model consisting of just three coupled oscillators, the authors show that such a description is sufficient to reproduce the wolf phenomenon. 
The study shows that the ratio between the vibration amplitudes at the left and the right bridge foot can be used as a reliable quantitative indicator to detect the wolf note in the time domain. 
Moreover, it demonstrates that the wolf can be effectively eliminated by actively increasing the body damping through an electromechanical feedback mechanism.

The most recent contribution among the works discussed here is \cite{gourc2022}, which provides a clear and up-to-date discussion of the physical origins of the wolf tone. Using a simplified bowed-string instrument model consisting of a linear string coupled to a single body resonance and excited by Coulomb friction, and without introducing any wolf-killing device, the authors focus exclusively on the emergence of the phenomenon. The governing equations are treated through modal decomposition rather than direct numerical time integration, allowing them to show that the wolf tone arises from a loss of stability of the Helmholtz motion through nonlinear bifurcations of periodic solutions.

Finally, let us mention the related approach \cite{mansour2017}, where the wolf note is indirectly investigated through playability criteria, such as the minimum bow force required to sustain Helmholtz motion, taking into account the measured body admittance of the instrument.


\paragraph{Main contribution.}
In this paper, we consider a mathematical model for the coupled dynamics consisting of 
\begin{itemize}
	\item a stiff string, including both tension-driven (second-order) and bending-stiffness (fourth-order) contributions. The string can be excited either by a plucking impulse or by a bowing force;
	\item a stiff body, described by a \emph{two-dimensional} thin plate, including again tension-driven (second-order) and bending-stiffness (fourth-order) contributions;  
	\item a bridge, modelled as a lumped mass with a single degree of freedom, elastically coupled to the string at one point and to the body at two distinct contact points representing the bridge feet. The two bridge-body connections are characterized by different stiffness values, accounting for the asymmetric loading induced by the string not acting at the bridge midline;
	\item one or two resonance suppressors, each of which is modelled as a mass-spring-damper system with a single degree of freedom, elastically and viscously coupled to the body at a specific point. Each suppressor can be placed anywhere on the body and it acts as a local resonator with a prescribed natural frequency, determined by its mass and stiffness parameters.
\end{itemize} 
Unlike most existing approaches, the proposed framework does not rely on modal decomposition: The coupled string-bridge-body-suppressor system is solved directly in the time domain by finite-difference discretization of the governing partial differential equations.
This fully numerical treatment allows us to handle both one- and two-dimensional subsystems within a unified framework, without prior modal identification or truncation. 
Moreover, our approach naturally extends to more general, in particular nonlinear, models, where 
(linear) modal analysis is no longer sufficient. Note that the approach operates downstream of the full system by relying solely on sampled vibrations of the instrument body. 


In order to investigate the effect of the suppressors, we consider a set of $\Nfreq$ possible lengths (or, equivalently, frequencies) of the string, corresponding to $N_f$ different notes that can be played through suitable finger placement on the string (The C string of a cello is taken as a reference).
Then, we use three performance indicators: 
the first one quantifies the wolf intensity, i.e.\ the degree of instability around the considered set of notes;
the second one measures the sound attenuation due to the wolf suppressor(s);  
finally, the third one measures the acoustic fidelity, defined as the deviation of the instrument's overall frequency response from that of the original noncorrected instrument. 
Guided by these three objectives, we determine the best placement and parameters of suppressor(s) that eliminate the wolf note while preserving as much as possible the tonal character of the instrument across the full spectrum.

\section{The model}\label{sec:model}
\begin{figure}[h!]
	\centering
	\includegraphics[width=0.5\textwidth]{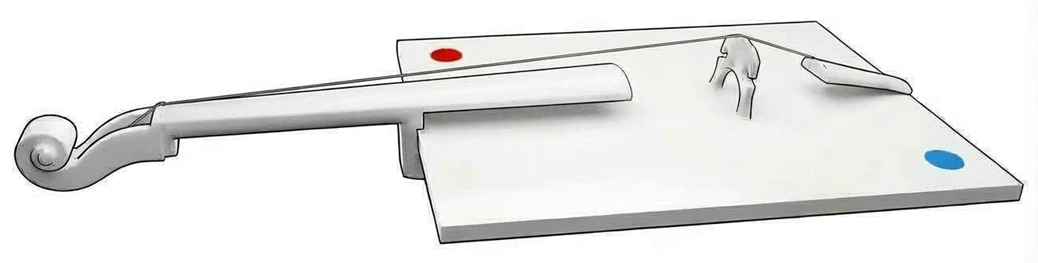}
	\caption{A pictorial representation of the string-bridge-body-suppressors system considered in the paper. 
	The red spot represents the registration point of the audio signal. 
	The blue one is the wolf suppressor.}
	\label{scheme_of_the_cello}
\end{figure}

In this section we present the model in detail. 
As already mentioned in the Introduction, the system is composed of a stiff string connected to a two-dimensional stiff body by means of a bridge represented by a single degree-of-freedom oscillator. 
The body, in turn, is connected to one or more wolf suppressors, see Figure \ref{scheme_of_the_cello}. 
The system is excited by plucking or bowing the string. 
For equations and their numerical approximation, we follow classical mathematical methods \cite{bilbaobook, bilbao2023}.

\subsection{Stiff string}
We consider a string of length $\ell$ [m] and we denote by $u(x,t)$ [m] its transverse displacement, with $u:[0,\ell]\times[0,+\infty]\to \mathbb R$.
We assume the string dynamics is given by
\begin{equation}\label{eq:stiffstring}
	u_{tt} = c_\st^2 u_{xx} 
	- r_\st^2 u_{xxxx} 
	-\beta_\st u_t 
	+ \frac{\Fstrbr}{\rho_\st A_\st}\delta(x-x_\br) 
	+ \frac{\Fattack}{\rho_\st A_\st}\delta(x-x_\att)\,,
\end{equation}
for $x\in[0,\ell]$ and $t>0$, where 
$\beta_\st$ [s$^{-1}$] is the velocity-proportional damping coefficient,
$\Fstrbr$ [N] is a point force accounting for the connection between the string and the bridge, 
$\Fattack$ [N] is another point force responsible for the string vibration, 
$x_\br$ is the connection point, on the string, with the bridge,
$x_\att$ is the point, on the string, excited with a pluck or the bow,
$\delta(x-x_0)$ is the Dirac delta centred in $x_0$, for any $x_0\in(0,\ell)$,
$\rho_\st$ [kg/m$^3$] is the material density,  
$A_\st$ [m$^2$] is the cross-sectional area,
$c_\st$ [m/s] is the wave speed, 
and $r^2_\st$ [m$^4$/s$^2$] is the bending coefficient. 
The last two fundamental quantities can be, in turn, derived by other quantities which are easier to measure, as   
\begin{equation}
c_\st = \sqrt{\frac{T_\st}{\rho_\st\, A_\st}} \qquad \text{and} \qquad  
r_\st = \sqrt{\frac{E_\st \, I_\st}{\rho_\st \, A_\st}}.
\end{equation}
where
$T_\st$ [N] is the string tension,  
$E_\st$ [Pa] is Young's modulus, and
$I_\st$ [m$^4$] is the second moment of area. 



\subsection{Stiff body}
We consider a plate of size $L$ $\times$ $L$ [m$^2$] and we denote by $w(x,y,t)$ [m] its transverse displacement, with $w:[0,L]^2\times[0,+\infty]\to \mathbb R$.
We assume the plate dynamics is given by
\begin{multline}\label{equation_body_w}
	w_{tt} = c_\bo^2 \nabla^2 w
	-r_\bo^2\nabla^4 w - 
	\beta_\bo w_t 
	+ \frac{\Fbobrleft}{\rho_\bo h_\bo}\delta(x-x_\brl,y-y_\brl) 
	+ \frac{\Fbobrright}{\rho_\bo h_\bo}\delta(x-x_\brr,y-y_\brr)\ + \\
	+ \sum_{\su} \frac{\Fbosupp}{\rho_\bo h_\bo}\delta(x-x_\su,y-y_\su)  \,,
\end{multline}
for $(x,y)\in[0,L]^2$ and $t>0$, 
where $\beta_\bo$ [s$^{-1}$] is the damping coefficient,
$\Fbobrleft$ and $\Fbobrright$ [N] are point forces exchanged between the body and the left and right foot of the bridge, respectively, 
$\Fbosupp$ [N] is another point force exchanged between the body and wolf suppressors (they can be many),
$\rho_\bo$ [kg/m$^3$] is the material density, 
$h_\bo$ [m] is the plate thickness,
$(x_\brl,y_\brl),(x_\brr,y_\brr)\in(0,L)^2$ 
are the connection points, on the plate, with the left and the right foot of the bridge, respectively,
and
$(x_\su,y_\su)\in(0,L)^2$ is the connection point, on the plate, with any of the suppressors,
$\delta(x-x_0,y-y_0)$ is the 2D Dirac delta centred in $(x_0,y_0)$, for any $(x_0,y_0)\in(0,L)^2$; 
finally, $\nabla^4w$ is a short notation for the bi-laplacian $w_{xxxx}+2w_{xxyy} + w_{yyyy}$. 

By also considering 
the Poisson's ratio $\nu$ [dimensionless], 
the Young's modulus $E_\bo$ [Pa], 
the (isotropic) in-plane membrane tension per unit length of the plate $T_\bo$ [N/m],
similarly to before we have
\begin{equation} 
c_\bo = \sqrt{\frac{T_\bo}{\rho_\bo\, h_\bo}} \qquad \text{and} \qquad 
r_\bo=\sqrt{\frac{D}{\rho_\bo h_\bo}}, \quad \text{with} \quad
D=\frac{E_\bo h^3_\bo}{12(1-\nu^2)}.
\end{equation}



\subsection{Bridge}
The bridge mediates the interaction between the string and the body. Denoting by $z_\textsc{br}(t)$ [m] its displacement from its equilibrium position and by $m_\textsc{br}$ [kg] its mass, the equation of motion is
\begin{equation}\label{eq:bridge}
	m_\textsc{br}\ \ddot z_\textsc{br}=
	- \Fstrbr
	- \Fbobrleft
	- \Fbobrright.
\end{equation}

\subsection{The wolf note suppressors}\label{sec:notch}
The wolf note suppressors are modelled as tuned mass dampers. Here we describe the dynamics of just one of them, but the full model will be able to deal with any number of suppressors, which are kept independent from each other. 
Denoting by $z_\su(t)$ [m] the displacement of the suppressor from its equilibrium position and by $m_\su$ [kg] its mass, the suppressor dynamics is driven by 
\begin{equation}
m_\su\ \ddot{z}_\su = -\Fbosupp. 
\end{equation}

\subsection{String-bridge connection}
The interaction force between the string and the bridge is given by
\begin{equation}
	\Fstrbr(t)
	:= k^\textsc{up}_\br\,\big(z_{\br}(t) - u(x_{\br},t)\big),
\end{equation}
where $k^\textsc{up}_\br$ [N/m] denotes the spring stiffness of the first spring in the bridge system.

\subsection{Body-bridge connection}
The interaction force between the plate and the bridge is obtained 
by summing the contributions of each connection point
\begin{equation}
	\Fbobrleft(t) := 
	k_\brl\,\big(z_{\br}(t) - w(x_\brl,y_\brl,t)\big),
	\qquad
	\Fbobrright(t) := 
	k_\brr\,\big(z_{\br}(t) - w(x_\brr,y_\brr,t)\big),
\end{equation}
where $k_\brl, k_\brr$ [N/m] denote the spring stiffness of the second and the third spring in the bridge system, respectively.

\subsection{Body-suppressor connection}
The interaction force between the body and all wolf suppressors is given by
\begin{equation}\label{def:Fbo-su}
\Fbosupp(t) := 
k_\su\big(z_\su - 
w(x_\su,y_\su,t)\big) 
+ 
\zeta_\su \big(\dot{z}_\su - w_t(x_\su,y_\su,t)\big)\,,
\end{equation} 
with $k_\su$ [N/m] and $\zeta_\su$ [kg/s] denoting the spring stiffness and the viscous damping coefficient of the suppressor, respectively.
It is useful to recall here that the following relationship holds true:
\begin{equation}\label{k=m2pif}
	k_\su=m_\su(2\pi f_\su)^2,
\end{equation} 
where $f_\su$ is the natural frequency of the suppressor.

\subsection{String excitation}
In our model, we assume that the string is initially at rest. To excite it, we apply a force $F_\att$ through one of the two standard techniques: plucking or bowing.

\subsubsection{Smooth plucking}
In the case of plucking, we define the smooth profile
\begin{equation}
	F_{\pl}(t) :=
	\left\{
	\begin{array}{ll}
		C_\pl \sin^2\left(\pi t/T_\pl\right), & t\leq T_\pl, \\
		0, & t>T_\pl,
	\end{array}
	\right.
\end{equation}  
where $C_\pl$ [N] is a positive constant and $T_\pl$ [s] is the duration of the pluck, typically very small.
Then, in \eqref{eq:stiffstring}, we set $F_\att = F_{\pl}$.

\subsubsection{Bowing}
In the case of bowing, the excitation force is much more complicated to describe since it is due to the continuous interaction between the bow hairs and the string. 
The bow hairs periodically grip the string for a short time and then release it, generating the well-known stick-and-slip dynamics.

Here, we use a model inspired by \cite{inacio2008}, mixed with the more classical approach described in \cite[Section 4.3.1]{bilbaobook}. 
At the initial time, we set the bow in stick mode (attack). 
At any time $t>0$, we compute the theoretical excitation force $F_\bow^*$ [N] needed to move the string at the exact speed of the bow, denoted by $V_\bow(t)$ (see Section \ref{sec:numerics} for details). 
The force $F_\bow^*$ is then compared with a threshold force $F_\bow^\textsc{max}$, which denotes the maximal tangential force the bow hairs can sustain, beyond which slipping occurs.
Denoting by $F_n(t)$ [N] the normal force applied by the bow and by
\begin{equation}\label{def:Vrel}
V_\rel(t):=u_t(x_\bow,t) - V_\bow(t) \qquad \text{[m/s]}
\end{equation}
the relative velocity between the string and the bow, we simply have 
\begin{equation}\label{def:Fbow}
    F_\bow=\left\{
    \begin{array}{ll}
    -F_n\ \mu_\stc \ \text{sign}(V_\rel), & |F_\bow^*|<F_\bow^\textsc{max}, \\ [1mm]
    -F_n\ \mu_\dyn \ \text{sign}(V_\rel), & \text{otherwise},
    \end{array}
    \right.
\end{equation}
where $\mu_\stc$ [dimensionless] and $\mu_\dyn$ [dimensionless] are the static and dynamic friction coefficients, respectively.
Note that the force $F_\bow^*$ is used to set the correct stick/slip phase only, while the force that actually describes the combined contribution of bow hairs is a Coulomb-type friction force. 
Finally, in \eqref{eq:stiffstring}, we set $F_\att = F_\bow$.

The advantage of the proposed approach is that it is quite stable from the numerical point of view, and does not require using the Newton method at every time step as suggested in, e.g., \cite{bilbaobook}.

\section{Wolf suppressor performance indicators}\label{sec:indicators}
Let us consider a discrete range of string lengths $\{\ell_1,\ldots,\ell_{N_f}\}$, for some $N_f>0$, which correspond to $N_f$ frequencies $\{f_1,\ldots,f_{N_f}\}$ as well as to $N_f$ musical notes that can be played on a single instrument string.
We also assume that one of those frequencies, say $f_\wolf\in\{f_1,\ldots,f_{N_f}\}$ generates the wolf note when no suppressors are employed.

Our goal will be to find optimal parameters and positions for the wolf suppressor(s) in such a way that 
\begin{itemize}
	\item the wolf tone is not audible on any musical note in the desired range;
	\item the wolf suppressor does not generate excessive attenuation of the response, i.e.\ a premature decay of the note. In fact, this is a well known drawback of some wolf suppressors caused by an anti-resonance effect which does not allow certain notes to be sustained well. Note that this effect \textit{adds} to that of damping, which is, of course, always present in real instruments;
	\item the spectra of all notes $f_i\neq f_\wolf$ in the range are not excessively modified, so as to keep the instrument as similar as possible to the original one (with no wolf suppressors);
\end{itemize}
To this end, we define three quantitative indicators to evaluate the system behavior: 
\begin{itemize}
	\item $J_\wolf$, for measuring the onset of the wolf note; 
    \item $J_\clos$, for measuring the attenuation of the waveform in the long term;
    \item $J_\fid$, for measuring the spectrum fidelity with the original instrument;
\end{itemize}
Note that the first and the third indicators rely on spectral analysis, while the second one is time-domain based.
\begin{remark}\label{exrem:f_e_f*}
    In the present framework, the audio signal is measured and recorded at some point \emph{on the body}, after the full coupling between the string, bridge, and body has taken place.
    It is plain that a string isolated (detached) from the instrument, with the same length $\ell_i$, would produce a sound with a fundamental frequency slightly different from $f_i$. 
\end{remark}
\begin{remark}\label{rem:choice_f*wolf}
	Following Equation \eqref{k=m2pif} and Remark \ref{exrem:f_e_f*}, $k_\su$ is chosen in such a way that 
	\begin{equation}\label{howtochooseksu}
	k_\su=m_\su(2\pi f_\wolf)^2.
	\end{equation} 
\end{remark}

\subsection{$J_\wolf$}\label{sec:Jwolf}
Let us fix a string frequency $f_i$ in the range of playable notes, a final time $T$ for the simulation, and choose a registration point on the plate $(x_\rec,y_\rec)\in(0,L)^2$.
We define
\begin{equation}\label{def:yi(t)}
	y_i(t):=w(x_\rec,y_\rec,t), \quad t\in[0,T], \quad i=1,\ldots,N_f,
\end{equation}
as the signal (waveform) recorded on the body, corresponding to the frequency $f_i$.
The complete pipeline used to get $J_\wolf$ is the following:
\begin{itemize}
	\item In order to remove trivial amplitude scaling effects, the recorded signal is first normalized in amplitude.
	We define
	\begin{equation}\label{ytilde_per_Jwolf}
		\widetilde y_i(t) =
		\frac{y_i(t)-\underbar y_i}{\widebar y_i -\underbar y_i}, \qquad
		\text{with} \quad
		\widebar y_i:=\max_{s\in[0,T]} y_i(s), \quad \text{and} \quad
		\underbar y_i:=\min_{s\in[0,T]} y_i(s).
	\end{equation}
	\item The analytic signal associated with $\widetilde y_i$ is defined as
	\begin{equation}
		z_i(t):=\widetilde y_i(t)+i\mathcal H\{\widetilde y_i\}(t),
	\end{equation}
	where $\mathcal H\{\cdot\}$ denotes the Hilbert transform.
	\item The instantaneous amplitude (envelope) of the response is defined as the modulus of the analytic signal,
	\begin{equation}
		a_i(t)=|z_i(t)|,\quad t\in[0,T].
	\end{equation}
    The envelope $a_i(t)$ captures amplitude modulations. 
	\item 
	To isolate the slow amplitude modulation component characteristic of the wolf-note beating, the envelope is filtered by convolution with a \emph{rectangular} kernel of width $\Theta$ (moving average filter),
	\begin{equation}
		K_\Theta(t)=\frac 1\Theta 
		\mathbf 1_{[-\Theta/2,\Theta/2]}(t).
	\end{equation}
	The filtered envelope is therefore
	\begin{equation}
		\overline a_i(t)= 
		(a_i\star K_\Theta)(t)
		=\frac 1\Theta \int_{t-\frac\Theta 2}^{t+\frac\Theta 2} 
		a_i(\tau)d\tau, \qquad t\in[0,T],
	\end{equation}
	where $a_i$ is extended outside $[0,T]$ by reflection at the endpoints. 
	\item The filtered envelope is detrended by subtracting its temporal mean,
	\begin{equation}
		e_i(t)=\overline a_i(t)-\langle \overline a_i\rangle,
		\qquad
		\langle \overline a_i\rangle=\frac1T\int_0^T \overline a_i(t) dt.
	\end{equation}
	\item We quantify the presence of low-frequency beating by measuring the fraction of spectral energy of $e_i$ concentrated in a prescribed band $[f_-,f_+]$.
	Denoting by $\mathcal F\{e_i\}$ the Fourier transform of $e_i$, we define
	\begin{equation}
		E^\band_i=
		\int_{f_-}^{f_+}\left|\mathcal F\{e_i\}(f)\right|^2 df,
		\qquad
		E^\tot_i = 
		\int_0^{f_\maxsc}\left|\mathcal F\{e_i\}(f)\right|^2 df,
	\end{equation}
	where $f_\maxsc$ is a prescribed upper bound used for normalization.
	\item The single wolf-note indicator associated with the $i$-th excitation frequency is defined as
	\begin{equation}\label{def:jiwolf}
		j_\wolf^i=\frac{E^\band_i}{E^\tot_i}.
	\end{equation}
	Values of $j_\wolf^i$ close to 1 correspond to a concentration of modulation energy in the low-frequency band $[f_-,f_+]$ and therefore to a pronounced wolf-note effect at frequency $f_i$.
	\item Finally, since we aim at preventing the onset of the wolf note on \emph{all} the considered frequencies, we define
	\begin{equation}\label{def:Jwolf}
		J_\wolf=\max_{i=1,\ldots,N_f} \{j^i_\wolf\}. \qquad [\%]
	\end{equation}
\end{itemize}

\subsection{$J_\clos$}\label{sec:Jclosing}
The second functional, named $J_\clos$, is introduced to quantify the residual amplitude of the waveform near the end of the simulation time. 
Its role is to measure possible anti-resonance effects induced by the wolf suppressor, which reduces beating at the cost of generating excessive attenuation of the response. 
By minimizing $J_\clos$, we discourage wolf suppressor placements that suppress the wolf tone through an unintended anti-resonant ``closure'' of the tone, thereby preserving its natural sustain.

We define $y_i(t)$ as in \eqref{def:yi(t)} and, as before, a final time $T$ for the simulation. Then, we define
\begin{equation}\label{def:Jclos}
	J_\clos=-\min_{i=1,\ldots,N_f}\max_{t\in[T_*,T]}\{|y_i(t)|\}, \quad [\text m]
\end{equation} 
where $T_*$ is a time reasonably close to the final one $T$. 
The term $\max_{t\in[T_*,T]}\{|y_i(t)|\}$ measures the degree of tone attenuation, as a small value of this term is necessarily associated with a ``closed'' tone.
Minimizing $J_\clos$ (i.e., maximizing $-J_\clos$), we are actually maximizing the minimal attenuation among all notes, thus avoiding excessive closures.
%
%
%
%
%
%
%
%

\subsection{$J_\fid$}\label{sec:Jfidelity}
In order to measure the difference between the overall acoustic response of the instrument with and without the wolf suppressors, we adopt the following pipeline starting from the recorded waveform $y_i(\cdot)$.
\begin{itemize}
	\item In order to remove trivial amplitude scaling effects, the recorded signal is first normalized in amplitude.
	We then define $\widetilde y_i$ as in \eqref{ytilde_per_Jwolf}.
	\item The magnitude spectrum is computed as 
	\begin{equation}
	A_i(f)=|\mathcal F\{\widetilde y_i\}(f)|, \quad f\geq 0. 
	\end{equation}
	\item To compress the dynamic range and improve the visibility of secondary spectral components, we introduce the logarithmic magnitude spectrum 
	\begin{equation}
	S_i(f)=20\log_{10}\big(\max\{A_i(f),\epsilon\}\big)\,,
	\end{equation}
	where $\epsilon>0$ is a small dimensionless parameter preventing the logarithmic singularity at zero. $S_i$ is taken as reference spectrum for the simulation without wolf suppressors.
	\item  When the wolf suppressors are employed, all previous computations are repeated, getting $S_i^\su(f)$. Then, the two spectra are compared using the standard $L^1$ norm
	\begin{equation}
	E_i=\int_0^{+\infty} |S_i(f)-S_i^\su(f)|df,
	\end{equation}
	and, finally, we compute the average of all frequencies different from $f_\wolf$
	\begin{equation}
	J_\fid=\frac{1}{N_f-1}\sum_{i\neq \wolf }E_i \quad [\text{dB$\cdot$Hz}].
	\end{equation}
\end{itemize}

%
%
%
%
%
%
%
%
\section{Numerical approximation}\label{sec:numerics}
Let us introduce structured and uniform numerical grids for the string and the body. 
Denote by $M_\st$ the number of cells for the string, and by $M_\bo \times M_\bo$ the number of cells for the body. 
We also denote by $N^T$ the number of time steps, and we define the time step $\Delta t=T/N^T$; the time step is the same for all subsystems (string, body, bridge, wolf suppressors).
Note that, while the number of time steps is defined \textit{a priori}, the space grid sizes $\Delta x_\st$, $\Delta x_\bo$, and $\Delta y_\bo$ are derived from $N^T$ and some model parameters, as to respect the CFL condition and minimizing numerical diffusion and dispersion.

As usual, we denote by $U_i^n$ the approximated value of the function $u$ at the grid node indexed by $i$, with $i=1,\ldots,M_\st$, at time step $n$, with $n=1,\ldots,N^T$; 
analogously, we denote by $W_{i,j}^n$ the approximated value of the function $w$ at the grid node indexed by $(i,j)$, with $i,j=1,\ldots,M_\bo$ at time step $n$, with $n=1,\ldots,N^T$.

\paragraph{String.}
For the string, we choose 
\begin{equation}
	\Delta x_\st = \sqrt{\frac{(c_\st \Delta t)^2 + \sqrt{(c_\st\Delta t)^4 + 16(r_\st\Delta t)^2}}{2}}\,,
\end{equation}
and we define
\begin{equation}
\lambda_\st = \left(\frac{c_\st \Delta t}{\Delta x_\st}\right)^2, \qquad
\mu_\st = \left(\frac{r_\st \Delta t}{\Delta x_\st^2}\right)^2, \qquad 
\tau_\st=\frac12 \beta_\st \Delta t. 
\end{equation}
The scheme for the string reads as 
\begin{multline}\label{schemestring}
	U_i^{n+1} = \frac{1}{(1+\tau_\st)} \bigg[
	(2 - 2\lambda_s - 6\mu_s) U_i^n + 
	\lambda_s (U_{i+1}^n + U_{i-1}^n) - 
	\mu_s (U_{i+2}^n - 4U_{i+1}^n - 4U_{i-1}^n + U_{i-2}^n) - 
	U_i^{n-1} + \\ 
	\tau_\st U_i^{n-1} +
	\frac{\FNattack \Delta t^2}{\rho_\st A_\st \Delta x_\st}\delta_{i,i_\att} +
	\frac{\FNstrbr \Delta t^2}{\rho_\st A_\st \Delta x_\st}\delta_{i,i_\br}\bigg]\,,
\end{multline} 
where $\delta_{\cdot,\cdot}$ is the Kronecker delta, and $i_\att$ and $i_\br$ are the indices corresponding to $x_\att$ and $x_\br$, respectively. 

We adopt simply supported boundary conditions, which are implemented as follows.
\begin{itemize}
	\item boundary nodes ($u=0$): 
	$
	U_2^n = 0, \quad 
	U_{M_\st-1}^n = 0,
	$
	\item ghost nodes ($u_{xx}=0$): 
	$
	U_1^n = -U_3^n, \quad
	U_{M_\st}^n = -U_{M_\st-2}^n.
	$
\end{itemize}


\paragraph{Body.} 
Regarding the body, we choose
\begin{equation}
	\Delta x_\bo = \Delta y_\bo = \max\left(\sqrt 2 c_\bo \Delta t \,, 2\sqrt{r_\bo \Delta t}\right),
\end{equation}
and we define
\begin{equation}
	\lambda_\bo = \left(\frac{c_\bo \Delta t}{\Delta x_\bo}\right)^2, \qquad
	\mu_\bo = \left(\frac{r_\bo \Delta t}{\Delta x_\bo^2}\right)^2, \qquad
	\tau_\bo=\frac12 \beta_\bo \Delta t.
\end{equation}
The scheme for the body reads as 
\begin{multline}
	W_{i,j}^{n+1} = \frac{1}{(1+\tau_\bo)} \bigg[
	2 W_{i,j}^n - W_{i,j}^{n-1} +
	\lambda_\bo (W_{i+1,j}^n - 2 W_{i,j}^n + W_{i-1,j}^n + W_{i,j+1}^n - 2 W_{i,j}^n + W_{i,j-1}^n)\ +\\
	-20 \mu_\bo W_{i,j}^{n} \ + \ 
	8 \mu_\bo ( W_{i+1,j}^{n} + W_{i-1,j}^{n} + W_{i,j+1}^{n} + W_{i,j-1}^{n} ) -
	2\mu_\bo (W_{i+1,j+1}^{n} + W_{i+1,j-1}^{n} + W_{i-1,j+1}^{n} + W_{i-1,j-1}^{n})\ + \\
	- \mu_\bo (W_{i+2,j}^{n} + W_{i-2,j}^{n} + W_{i,j+2}^{n} + W_{i,j-2}^{n}) +
	\tau_\bo W_{i,j}^{n-1} +
	\frac{\FNbobrleft \Delta t^2}{\rho_\bo h_\bo \Delta x^2_\bo}\delta_{i,i_\brl} \delta_{j,j_\brl} + \\
	\frac{\FNbobrright \Delta t^2}{\rho_\bo h_\bo \Delta x^2_\bo}\delta_{i,i_\brr} \delta_{j,j_\brr} +
	\sum\frac{\FNbosupp \Delta t^2}{\rho_\bo h_\bo \Delta x^2_\bo}\delta_{i,i_\su} \delta_{j,j_\su}
	\bigg]\,,
\end{multline} 
where $(i_\brl,j_\brl)$, $(i_\brr,j_\brr)$, and $(i_\su,j_\su)$ are the indices corresponding to $(x_\brl,y_\brl)$, $(x_\brr,y_\brr)$, and $(x_\su,y_\su)$, respectively.

Boundary conditions are treated as in the string case.


\paragraph{Bridge.} 
The scheme for the bridge reads as
\begin{equation}
	Z_\br^{n+1} = 2Z_\br^n - Z_\br^{n-1} - \frac{\Delta t^2}{m_\br} 
	(\FNstrbr+\FNbobrleft+\FNbobrright).
\end{equation}


\paragraph{Wolf suppressors.}
The scheme for each wolf suppressor reads as
\begin{equation}
	Z_\su^{n+1} = 2Z_\su^n - Z_\su^{n-1} - \frac{\Delta t^2}{m_\su} 
	\FNbosupp.
\end{equation}
$\dot z_\su$ in \eqref{def:Fbo-su} is approximated by 
$$
\frac{Z_\su^{n+1}-Z_\su^{n-1}}{2\Delta t}.
$$

\paragraph{Bowing.}
For the computation of the excitation force in the case of bowing, three remarks are in order: 
\begin{itemize}
    \item The value $F_\bow^*$ is computed directly at a discrete level using \eqref{schemestring}. We set $F_\bow^*$ equal to the value of $\FNattack$ in \eqref{schemestring} such that $(U^{n+1}_i-U^n_i)/\Delta t=V_\bow$. 
    \item To enhance numerical stability, we relax the sign function in \eqref{def:Fbow} defining
    $$
    \text{sign}(x)=\left\{ \begin{array}{ll} 1, & x>\varepsilon, \\ 0, & |x|<\varepsilon, \\ -1, & x<-\varepsilon. \end{array} \right.
    $$
    \item Both to enhance numerical stability and to be more faithful to real bows and cellos, we assume that the contact between the bow and the string happens on a surface with a certain extension. More precisely, we spread the bow force on three numerical nodes using a triangular kernel with weights $(1/4,1/2,1/4)$.  
\end{itemize}

Finally, the performance indicators introduced in Section \ref{sec:indicators} are approximated by standard numerical methods, FFT and inverse FFT. 

\clearpage
\section{Numerical tests}\label{sec:numericaltests}
In this section, we present some numerical tests for the full system excited either by pluck or bow, with 0, 1, and 2 suppressors. 
The goal is to study the effect of the suppressors as a function of their position on the body, observing the behaviour of the waveforms, the spectra, and the three indicators introduced in Section \ref{sec:indicators}.

The values of all (physical and numerical) parameters are summarized in Table \ref{tab:parameters}.
Lengths and frequencies introduced in Section \ref{sec:indicators} are summarized in Table \ref{tab:frequencies}.
Note that, while the tested frequencies correspond to real-life notes, the string lengths must differ from the physically measured lengths on a real instrument, owing to the simplifying assumptions underlying the model.
Also, the need for slightly different effective string lengths under plucked and bowed excitation is interpreted as a modeling artifact arising from the simplified stick-slip interaction and from the different nonlinear dynamical equilibria established in the two excitation regimes.
In Section \ref{sec:extrawolfs}, one additional frequency will also be considered. 
\begin{table}[h!]
	\caption{Parameters used in simulations.}
	\centering
	\begin{subtable}[t]{0.32\textwidth}
		\centering \vskip-2.1cm
		\begin{tabular}{|l|l|l|}
			\hline
			symbol        & value                & unit     \\ \hline\hline 
			$T_\st$	      & $1.2\times 10^2$     & N        \\ \hline
			$\rho_\st$	  & $7.8\times 10^3$     & kg/m$^3$ \\ \hline
			$A_\st$	      & $1.8\times 10^{-6}$  & m$^2$    \\ \hline
			$E_\st$	      & $2.0\times 10^{11}$  & Pa       \\ \hline
			$I_\st$	      & $9.8\times 10^{-14}$ & m$^4$    \\ \hline
			$\beta_\st$   & 0                    & s$^{-1}$ \\ \hline
		\end{tabular}
		\caption{String}
	\end{subtable}
	\begin{subtable}[t]{0.32\textwidth}
		\centering 
		\begin{tabular}{|l|l|l|}
			\hline
			symbol        & value                & unit     \\ \hline\hline  
			$L$           & $5.0\times 10^{-1}$  & m        \\ \hline 
                $T_\bo$	      & $7.5\times 10^4$     & N/m      \\ \hline
			$\rho_\bo$	  & $4.7\times 10^{2}$   & kg/m$^3$ \\ \hline
			$E_\bo$	      & $1.0\times 10^{10}$  & Pa       \\ \hline
			$h_\bo$	      & $4.0\times 10^{-3}$  & m  \\ \hline
			$\nu$	      & $2.5\times 10^{-1}$  & --       \\ \hline
                $\beta_\bo$   & 0                    & s$^{-1}$ \\ \hline
		\end{tabular} 
		\caption{Body}
	\end{subtable}
	\begin{subtable}[t]{0.32\textwidth}
		\centering \vskip-2.1cm
		\begin{tabular}{|l|l|l|}
			\hline
			symbol        & value                & unit     \\ \hline\hline  
			$m_\br$	      & $2.0\times 10^{-2}$  & kg       \\ \hline
			$m_\su$	      & $8.5\times 10^{-3}$  & kg       \\ \hline
			$k^\textsc{up}_\br$	   & $4.9\times 10^2$ & N/m \\ \hline
			$k^\textsc{left}_\br$  & $7.0\times 10^4$ & N/m \\ \hline
			$k^\textsc{right}_\br$ & $3.0\times 10^4$ & N/m \\ \hline
			$\zeta_\su$	  & 2.1                  & kg/s     \\ \hline
		\end{tabular}
		\caption{Bridge and wolf suppressors}
		\label{tab:parameters_bridge_notch}
	\end{subtable}
	\begin{subtable}[t]{0.32\textwidth}
		\centering
		\begin{tabular}{|l|l|l|}
			\hline
			symbol        & value                & unit \\ \hline\hline  
			$C_\pl$	      & 1                    & N    \\ \hline
			$T_\pl$	      & $4.55\times 10^{-3}$ & s    \\ \hline
			$V_\bow$	  & $2.0\times 10^{-1}$  & m/s  \\ \hline
			$F_n$	      & $1.0$       & N    \\ \hline
			$F_\bow^\textsc{max}$	  & $2.5$    & N  \\ \hline
			$\mu_\stc$	  & $6.0\times 10^{-1}$  & --   \\ \hline
			$\mu_\dyn$	  & $2.0\times 10^{-1}$  & --   \\ \hline
            $\varepsilon$ & $1.0\times 10^{-2}$  & m/s   \\ \hline
		\end{tabular}
		\caption{Pluck and bow}
	\end{subtable}
	\begin{subtable}[t]{0.32\textwidth}
		\centering \vskip-2.35cm
		\begin{tabular}{|l|l|l|}
			\hline
			symbol     & value               & unit \\ \hline\hline  
			$\Delta t$ & $5.7\times 10^{-6}$ & s    \\  \hline
			$T$        & $1.0$               & s   \\ \hline
			$T_*$      & $0.9$               & s   \\ \hline
			$\Theta$   & $1.0\times 10^{-2}$ & s    \\ \hline
			$f_-$      & $2.0$               & Hz   \\ \hline
			$f_+$      & $1.3\times 10$      & Hz   \\ \hline
			$f_\maxsc$ & $1.0\times 10^{2}$  & Hz   \\ \hline
		\end{tabular}
		\caption{Numerics and thresholds}
	\end{subtable}
	\begin{subtable}[t]{0.32\textwidth}
		\centering  \vskip-2.35cm
		\begin{tabular}{|l|l|l|}
			\hline
			point             & value    & unit \\ \hline\hline  
			$x_\att$          & 50       & --   \\  \hline
			$x_\br$           & 70       & --   \\  \hline
			$(x_\brl,y_\brl)$ & (42,48)  & --   \\ \hline
			$(x_\brr,y_\brr)$ & (42,52)  & --   \\ \hline
			$(x_\su,y_\su)$   & variable & --   \\ \hline
			$(x_\rec,y_\rec)$ & (42,18)  & --   \\ \hline
		\end{tabular}
		\caption{Contact points (\% on maximal size)}
		\label{tab:parameters_contactpoints}
	\end{subtable}
	\label{tab:parameters}
\end{table}
\begin{table}[h!]
	\caption{Lengths $\{\ell_i\}_i $ of the string, and corresponding frequencies $\{f_i\}_i$ and notes as sampled on body. The 5-th one generates the wolf tone.}
	\centering
	\begin{tabular}{|l|c|c|c|c|c|c|c|c|c|}
		\hline
		& 1 & 2 & 3 & 4 & \textbf 5 & 6 & 7 & 8 & 9  \\ \hline  \hline 
            $\ell$ [cm] (pluck) & 24.8 & 23.4 & 22.2 & 20.9 & \textbf{19.7} & 18.9 & 17.8 & 16.9 & 16.0   \\ \hline
            $\ell$ [cm] (bow) & 25.1 & 23.8 & 22.6 & 21.6 & \textbf{20.1} & 18.6 & 17.7 & 16.8 & 15.9   \\ \hline\hline  
		$f$ [Hz]  & 196.0 & 207.2 & 220.0 & 233.1 & \textbf{246.9} & 261.6 & 277.2 & 293.7 & 311.1  \\ \hline
		note & G3 & G$\sharp 3$ & A3 & A$\sharp 3$ & \textbf{B3} & C4 & C$\sharp 4$ & D4 & D$\sharp 4$  \\ \hline
	\end{tabular}
	\label{tab:frequencies}
\end{table}

%
%
%
%
%
\subsection{Excitation by pluck}
In this section, we assume that the string is excited by a pluck and that the system is not damped. 
This is done to better highlight the effect of wolf suppressor(s).

\subsubsection{Test PLUCK-0S: No wolf suppressors}
Here we show how the model described in Section \ref{sec:model} is able to reproduce sounds for the 9 notes reported in Table \ref{tab:frequencies}, with no wolf suppressors. 
Among the considered notes, it turns out that frequency $f_5$ provides the wolf tone. 
We also show how the indicator $J_\wolf$ is able to capture the wolf tone among all the notes.

Note that the spectra computed here with no wolf suppressors will be used as a reference for computing $J_\fid$ in the following sections, where the wolf suppressor(s) will be enabled. 

\begin{figure}[h!]
	\centering
	\begin{subfigure}{0.49\textwidth}
		\centering
		\includegraphics[width=\textwidth]{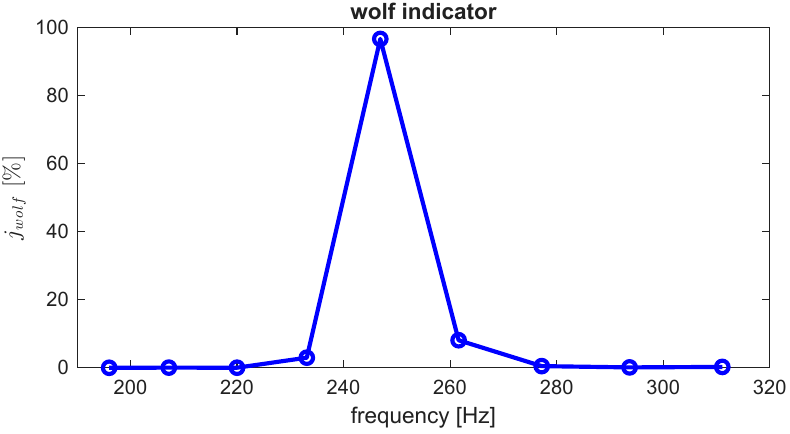} 
		\caption{Indicator $f_i\to j^i_\wolf$}
		\label{fig:pluck_0notch_allf_Jwolf}
	\end{subfigure}
	\\
	\begin{subfigure}{0.49\textwidth}
		\centering
		\includegraphics[width=\textwidth]{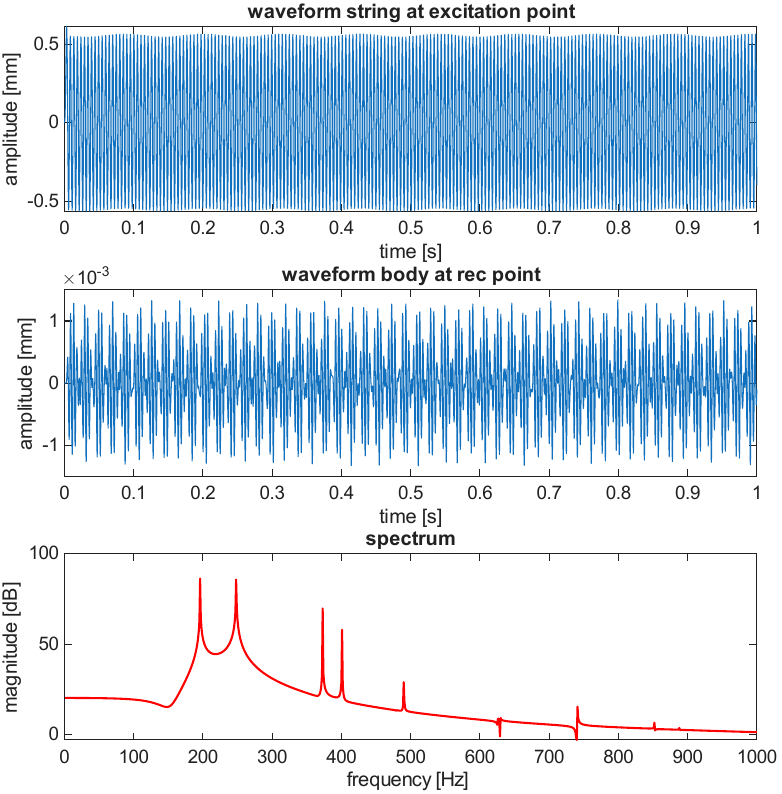}
		\caption{
			Waveform and spectrum for $f_1$ 	
                                   \audionotelocal{./supplementary_material/audio/9n_pluck_0notch_f1_audio.wav}
            \audionoteonline{www.emilianocristiani.it/attach/paper_wolfnote/audio/9n_pluck_0notch_f1_audio.wav}
		}
		\label{fig:pluck_0notch_f1_WF_SP}
	\end{subfigure}
	\begin{subfigure}{0.49\textwidth}
		\centering
		\includegraphics[width=\textwidth]{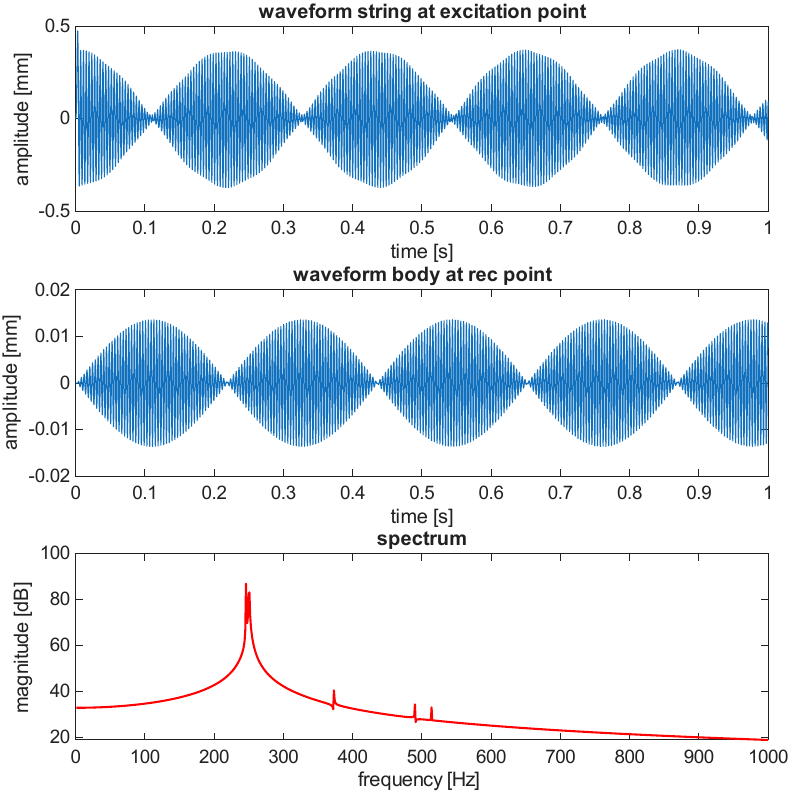}
		\caption{
			Waveform and spectrum for $f_5$
                                       \audionotelocal{./supplementary_material/audio/9n_pluck_0notch_f5_audio.wav} 
			\audionoteonline{www.emilianocristiani.it/attach/paper_wolfnote/audio/9n_pluck_0notch_f5_audio.wav}
		}
		\label{fig:pluck_0notch_f5_WF_SP}
	\end{subfigure}
	\caption{
		PLUCK-0S: Indicator $f_i\to j^i_\wolf$ and results for two frequencies from Table \ref{tab:frequencies}. 
		Waveforms of all frequencies are available 
                                         \href{./supplementary_material/large_figures/9n_pluck_0notch_allf_WF.pdf}{here-local} and
		\href{http://www.emilianocristiani.it/attach/paper_wolfnote/large_figures/9n_pluck_0notch_allf_WF.pdf}{here-online}.
		Sounds of all frequencies are available 
		                               \href{./supplementary_material/audio/9n_pluck_0notch_allf_audio.wav}{here-local} and
		\href{http://www.emilianocristiani.it/attach/paper_wolfnote/audio/9n_pluck_0notch_allf_audio.wav}{here-online}.
	}  
	\label{fig:pluck_0notch}
\end{figure}

Figure \ref{fig:pluck_0notch_allf_Jwolf} shows the indicator $j^i_\wolf$, defined in \eqref{def:jiwolf}, for any frequency $i=1,\ldots,N_f$. We see that $i^5_\wolf>95\%$, meaning that the wolf note is captured with very high precision. 

Figures \ref{fig:pluck_0notch_f1_WF_SP}-\ref{fig:pluck_0notch_f5_WF_SP} show waveforms and spectra of $f_1$ and $f_5$, respectively. 
Waveforms are recorded at the excitation point on the string and at the registration point on the body.
For $f_5$, both the wolf tone in the waveforms and the double peak in the spectrum are clearly visible. 
Note that the wolf-note beating is observable not only in the body waveform, but also in the string waveform, in response of the vibro-acoustic coupling between string and body. 
%
\clearpage
\subsubsection{Test PLUCK-1S: optimal placement of one wolf suppressor}\label{sec:testpluck-1S}
Here we enable one wolf suppressor and investigate to what extent it is able to prevent the onset of wolf note. 
We also study how the suppressor perturbs notes originally not affected by the wolf.
To begin with, we fix suppressor parameters $m_\su$ and $\zeta_\su$ as specified in Table \ref{tab:parameters}, $f_\su=f_\wolf=246.9$ (see Remark \ref{rem:choice_f*wolf}), and vary the suppressor position along the whole instrument body. More precisely, we test $45 \times 45$ suppressor positions, corresponding to the number of spatial nodes in the body grid.

\begin{figure}[h!]
	\centering
	\begin{subfigure}{0.33\textwidth}
		\centering
		\includegraphics[width=\textwidth]{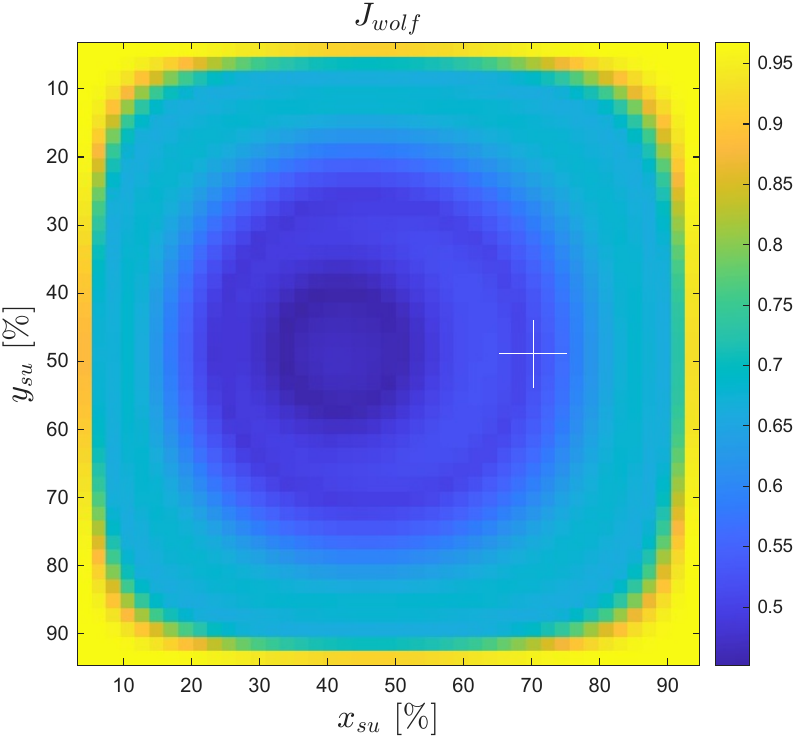} 
		\caption{$J_\wolf$}
		\label{fig:pluck_1notch_allf_Jwolf_landscape}
	\end{subfigure}
	\begin{subfigure}{0.33\textwidth}
		\centering
		\includegraphics[width=1.01\textwidth]{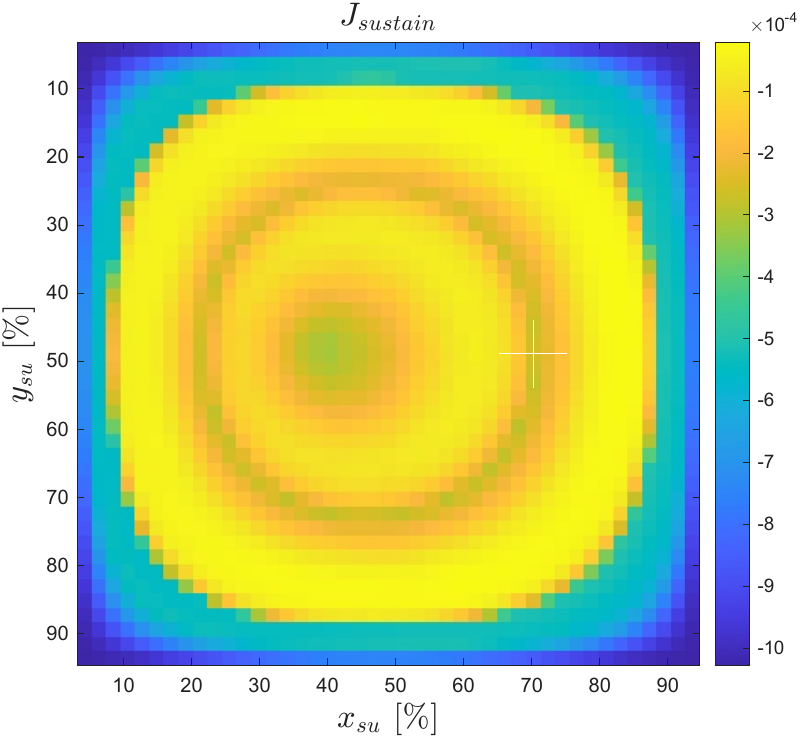} 
		\caption{$J_\clos$}
		\label{fig:pluck_1notch_allf_Jclos_landscape}
	\end{subfigure}
    \begin{subfigure}{0.33\textwidth}
		\centering
		\includegraphics[width=0.98\textwidth]{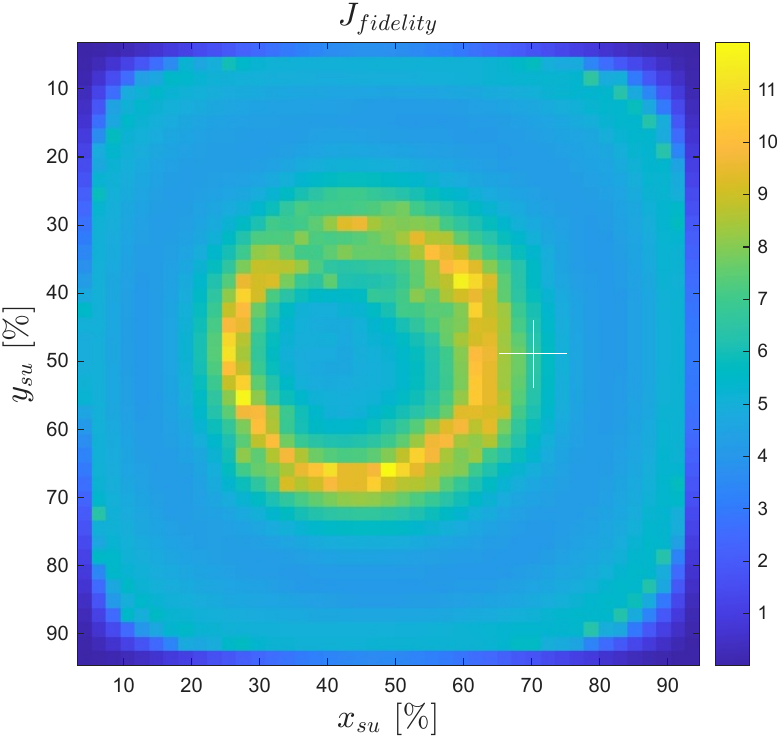} 
		\caption{$J_\fid$}
		\label{fig:pluck_1notch_allf_Jfid_landscape}
	\end{subfigure}
	\caption{PLUCK-1S. Heat maps of $J_\wolf$, $J_\clos$, and $J_\fid$ as a function of the wolf suppressor position $(x_\su,y_\su)$ on the instrument body (expressed as $\%$ of the body size). The white cross indicates the point chosen as best placement of the suppressor.
	}
	\label{fig:pluck_1notch_allf_allJ_landscape}
\end{figure}

Figure \ref{fig:pluck_1notch_allf_allJ_landscape} shows the heat maps of the indicators $J_\wolf$, $J_\clos$, and $J_\fid$ as a function of the suppressor position $(x_\su,y_\su)$.  

We immediately observe that the images have no radial symmetry with respect to the center of the body. The symmetry breaking is due to the bridge position, which is not located in the center, and to the bridge feet, which are not equal, cfr.\ Tables \ref{tab:parameters_bridge_notch} and \ref{tab:parameters_contactpoints}. 
We also see that the behaviour of the three indicators is different and, more importantly, they reach the minimum values at different points. 
In general, we can say that if we try to improve one indicator, we will likely worsen the others.
Also, we think that summing up the three indicators in order to get a single global indicator is not advisable because, in real life, the weight to give to each term depends on the piece currently being performed and it is a very subjective choice of the musician.\footnote{Probably the order of importance is $J_\wolf>J_\clos>J_\fid$ but the values to assign to the three possible weights $\alpha_1,\alpha_2,\alpha_3$ in order to define a global indicator $J=\alpha_1 J_\wolf+\alpha_2 J_\fid+\alpha_3 J_\clos$ are largely arbitrary.}

For the following tests, we heuristically choose the point $(70\%,49\%)$ (indicated with the white cross in the figure) as the best compromise. Indeed, a wolf suppressor located there is able to eliminate the wolf tone while still preserving an acceptable sustain and, thirdly in importance, fidelity. 

\begin{figure}[h!]
	\centering
	\includegraphics[width=\textwidth]{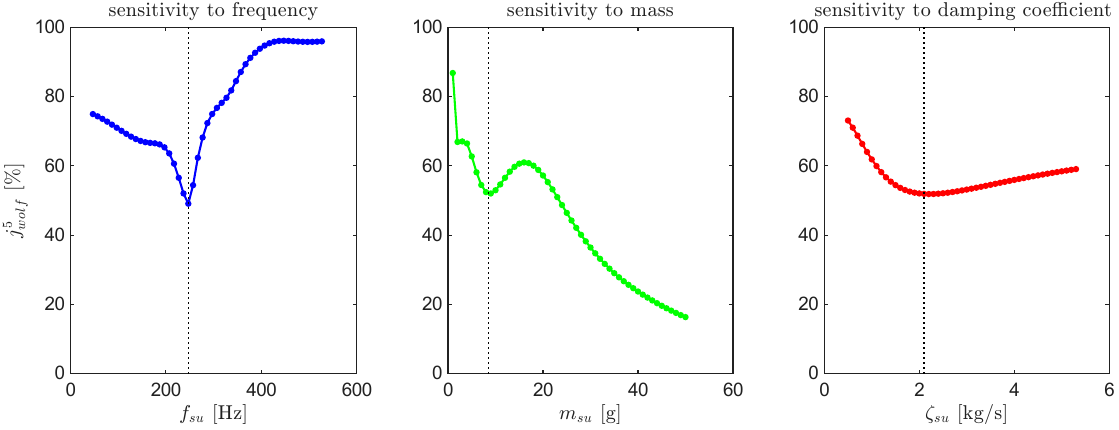} 
	\caption{PLUCK-1S. $j_\wolf^5$ as a function of suppressor parameters $f_\su$, $m_\su$, and $\zeta_\su$, respectively. 
	}
	\label{fig:pluck_1notch_f5_suppressorparameterssensitivy}
\end{figure}

Figure \ref{fig:pluck_1notch_f5_suppressorparameterssensitivy} shows the value of the indicator $j_\wolf^5$, corresponding specifically to the wolf note $f_5$, as a function of suppressor parameters $f_\su$, $m_\su$, and $\zeta_\su$, respectively. 
Each plot is obtained by varying one parameter at a time while fixing the other two.
Interestingly, one can observe that the suppressor parameters were chosen in such a way that each of them allows reaching the minimum value of $j_\wolf^5$ if the other parameters are kept fixed.
Considering that the optimal suppressor placement depends on the suppressor parameters and, in turn, the optimal choice of the suppressor parameters depends on the suppressor placement, we conclude that the considered 5-tuple $(x_\su,y_\su,f_\su,m_\su,\zeta_\su)$ represents a Nash-like locally optimal configuration with respect to single-parameter variations. In other words, any modification of \textit{a single} value leads to a higher value of the wolf indicator. 
On the other hand, it is possible that modifying two or more parameters at the same time allows for further improvement of the system performance.

\begin{figure}[h!]
	\centering
	\begin{subfigure}{0.49\textwidth}
		\centering
		\includegraphics[width=\textwidth]{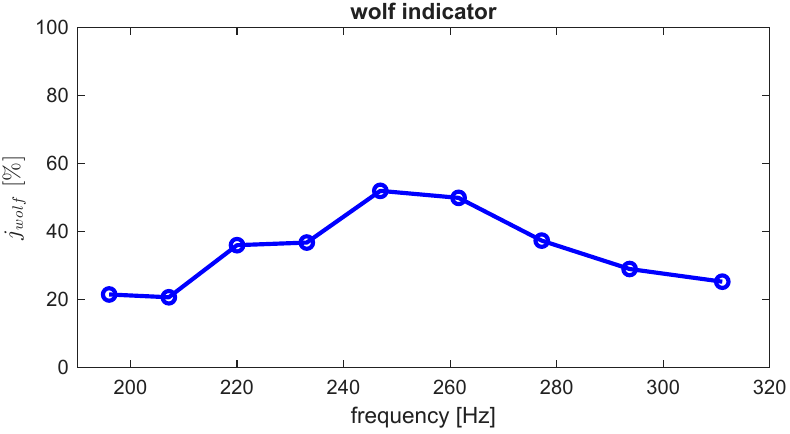} 
		\caption{Indicator $f_i\to j^i_\wolf$}
		\label{fig:pluck_1notch3323_allf_Jwolf}
	\end{subfigure}
	\\
	\begin{subfigure}{0.49\textwidth}
		\centering
		\includegraphics[width=\textwidth]{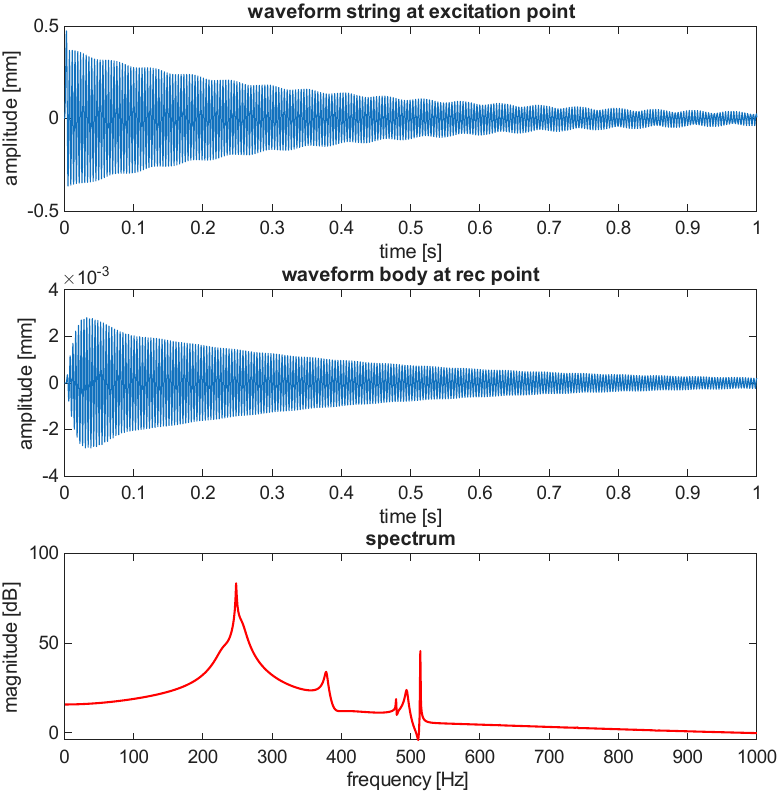}
		\caption{
			Waveform and spectrum for $f_5$, good placement of the wolf suppressor
			                        \audionotelocal{./supplementary_material/audio/9n_pluck_1notch3323_f5_audio.wav}
			\audionoteonline{www.emilianocristiani.it/attach/paper_wolfnote/audio/9n_pluck_1notch3323_f5_audio.wav}
		}
		\label{fig:pluck_1notch3323_f5_WF_SP}
	\end{subfigure}
	\begin{subfigure}{0.49\textwidth}
		\centering
		\includegraphics[width=\textwidth]{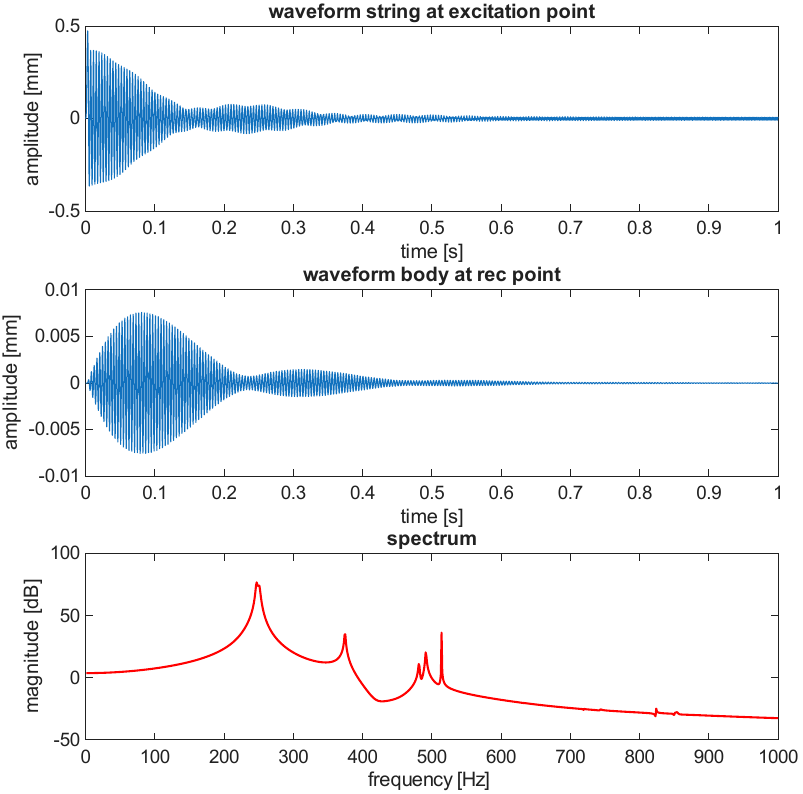}
		\caption{
			Waveform and spectrum for $f_5$, bad placement of the wolf suppressor
			                        \audionotelocal{./supplementary_material/audio/9n_pluck_1notch3610_f5_audio.wav} 
			\audionoteonline{www.emilianocristiani.it/attach/paper_wolfnote/audio/9n_pluck_1notch3610_f5_audio.wav}
		}
		\label{fig:pluck_1notch3610_f5_WF_SP}
	\end{subfigure}
	\caption{
		PLUCK-1S: Indicator $f_i\to j^i_\wolf$ and results for two frequencies from Table \ref{tab:frequencies}. 
		Waveforms of all frequencies are available 
		                            \href{./supplementary_material/large_figures/9n_pluck_1notch3323_allf_WF.pdf}{here-local} and
		\href{http://www.emilianocristiani.it/attach/paper_wolfnote/large_figures/9n_pluck_1notch3323_allf_WF.pdf}{here-online}.
		Sounds of all frequencies are available 
		                               \href{./supplementary_material/audio/9n_pluck_1notch3323_allf_audio.wav}{here-local} and
		\href{http://www.emilianocristiani.it/attach/paper_wolfnote/audio/9n_pluck_1notch3323_allf_audio.wav}{here-online}.
	}  
	\label{fig:pluck_1notch}
\end{figure}

Figure \ref{fig:pluck_1notch3323_allf_Jwolf} shows the indicator $j^i_\wolf$ for any frequency $i=1,\ldots,N_f$. We see that all frequencies are below 60\%, meaning that the wolf note is well suppressed. 

Figure \ref{fig:pluck_1notch3323_f5_WF_SP} shows waveforms and spectrum of $f_5$. 
We see that the original wolf tone is now completely suppressed, although the tone is attenuated. 
Nevertheless, these suppressor performances are acceptable in comparison with what we can get with a \textit{wrong} placement of the suppressor itself. 
For example, choosing a point in the yellow region in Figure \ref{fig:pluck_1notch_allf_Jclos_landscape}), we get the waveform shown in Figure \ref{fig:pluck_1notch3610_f5_WF_SP}. 
Here we have both the wolf note and a very strong attenuation.
%

\clearpage 
\subsection{Excitation by pluck for one extra wolf note}\label{sec:extrawolfs}
\begin{figure}[h!]
	\centering
	\includegraphics[width=\textwidth]{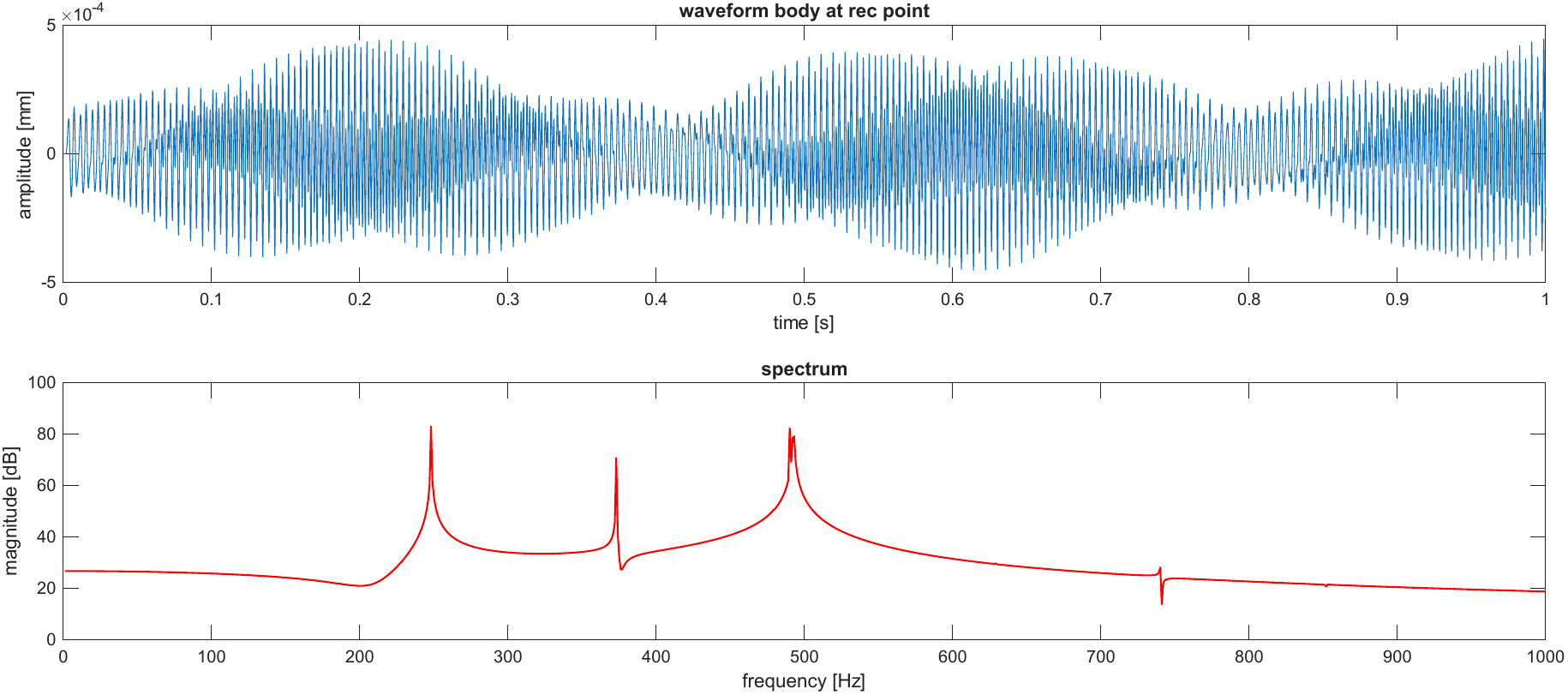}
	\caption{Waveform and spectrum for $f=492$ Hz.
                                       \audionotelocal{./supplementary_material/audio/9n_pluck_0notch_f492_audio.wav} 
			\audionoteonline{www.emilianocristiani.it/attach/paper_wolfnote/audio/9n_pluck_0notch_f492_audio.wav}
            }
	\label{fig:pluck_0notch_f492_WF_SP}
\end{figure}

For extra fun, we run the simulation for one additional string length, namely $\ell=10.4$ cm ($f=492$ Hz), which also produce the wolf note. 
This note falls outside the interval of notes considered up to now and that are summarized in Table \ref{tab:frequencies}.
Figure \ref{fig:pluck_0notch_f492_WF_SP} shows the waveform recorded on body and the spectrum for the extra frequency. 

\begin{figure}[h!]
	\centering
	\begin{subfigure}{0.45\textwidth}
		\centering
		\includegraphics[width=\textwidth]{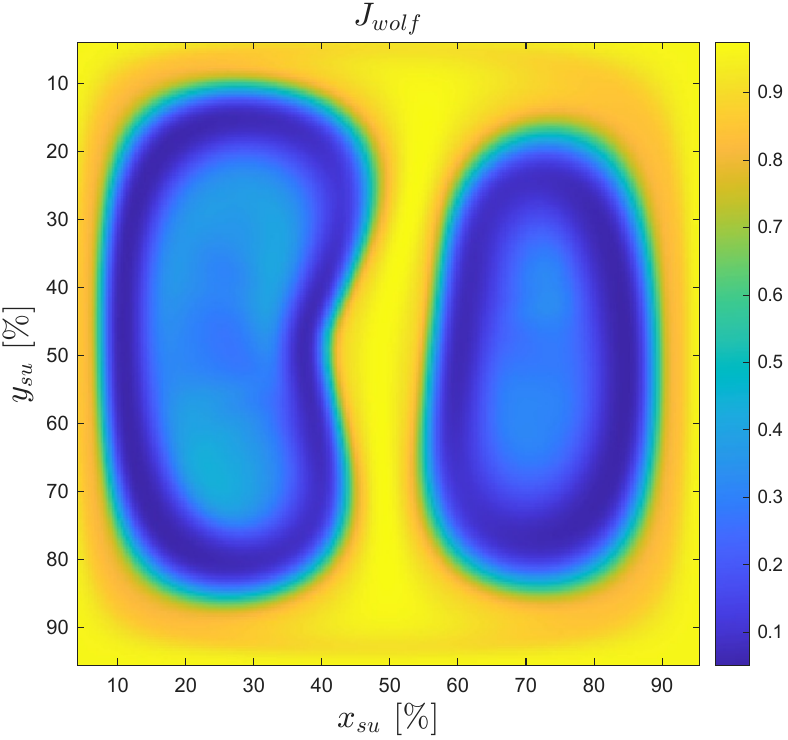} 
		\caption{$J_\wolf$ with $f_\wolf=246.9$ Hz}
		\label{fig:pluck_1notch_f492_Jwolf_landscape}
	\end{subfigure}
	\begin{subfigure}{0.45\textwidth}
		\centering
		\includegraphics[width=\textwidth]{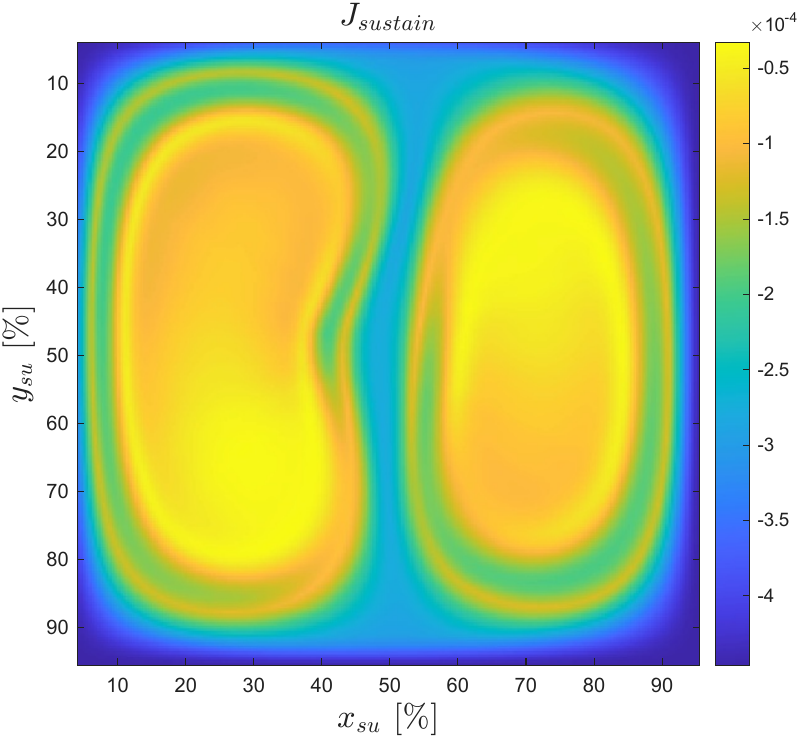} 
		\caption{$J_\clos$ with $f_\wolf=246.9$ Hz}
		\label{fig:pluck_1notch_f127_Jclos_landscape}
	\end{subfigure}
	\begin{subfigure}{0.45\textwidth}
		\centering
		\includegraphics[width=\textwidth]{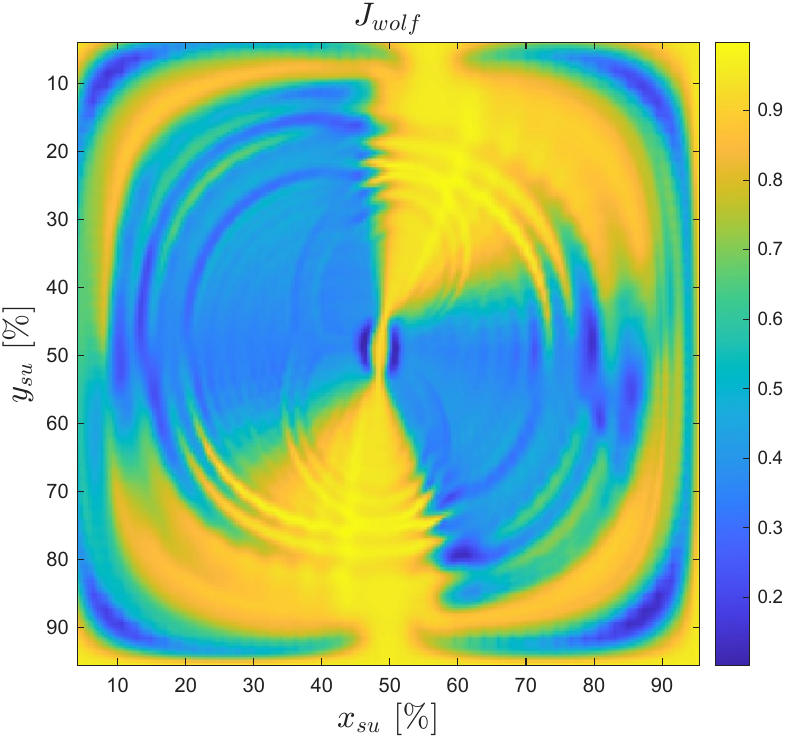} 
		\caption{$J_\wolf$ with $f_\wolf=492$ Hz}
		\label{fig:pluck_1notch_f492_n492_Jwolf_landscape}
	\end{subfigure}
	\begin{subfigure}{0.45\textwidth}
		\centering
		\includegraphics[width=\textwidth]{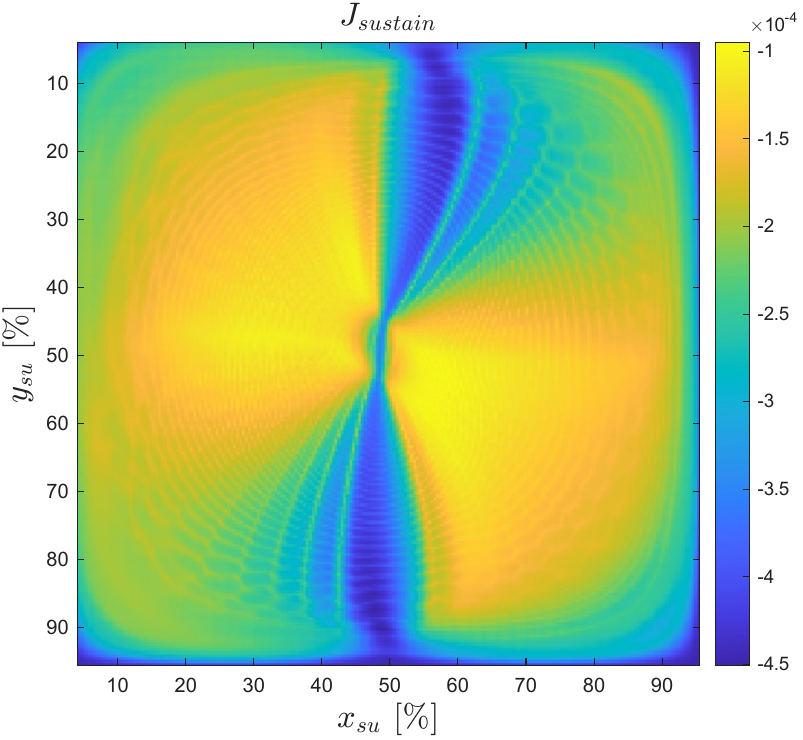} 
		\caption{$J_\clos$ with $f_\wolf=492$ Hz}
		\label{fig:pluck_1notch_f492_n492_Jclos_landscape}
	\end{subfigure}
	\caption{Heat maps of $J_\wolf$ and $J_\clos$ for $f=492$ Hz, as a function of the wolf suppressor position $(x_\su,y_\su)$ on the instrument body (expressed as $\%$ of the body size). $k_\su$ in \eqref{howtochooseksu} is chosen assuming either $f_\wolf=246.9$ Hz or $f_\wolf=492$ Hz. 
	}
	\label{fig:pluck_1notch_f492_allJ_landscape}
\end{figure}

Figure \ref{fig:pluck_1notch_f492_allJ_landscape} shows the heat maps of the indicators $J_\wolf$ and $J_\clos$ as a function of the suppressor position $(x_\su,y_\su)$ for the extra frequency (here $j_\wolf=J_\wolf$). 
We tested two distinct suppressor's resonance frequencies in \eqref{howtochooseksu}, namely $f_\wolf=246.9$ Hz and $f_\wolf=492$ Hz.  
%
The resulting structures are different from those in Figure \ref{fig:pluck_1notch_allf_allJ_landscape}, suggesting that the beating observed at 492 Hz does not involve the first body mode, but rather arises from the interaction with a higher body mode.
Note that for these computations we tested $215 \times 215$ suppressor positions. 
Since this number is larger than that of the body grid nodes, we use bilinear interpolation to distribute $\Fbosupp$ among the four grid points closest to each tested suppressor position.

\clearpage
\subsection{Excitation by bow}
In this section, we assume that the string is excited by a bow. 

\subsubsection{Test BOW-0S: No wolf suppressor}\label{sec:testbow-0S}
Here we show how the model described in Section \ref{sec:model} is able to reproduce sounds for the 9 notes reported in Table \ref{tab:frequencies}, with no wolf suppressors. 
Among the considered notes, the frequency $f_5$ provides the wolf tone. 
We also show how the indicator $J_\wolf$ is able to capture the wolf tone among all the notes.

As we have done for previous tests, the spectra computed here with no wolf suppressors will be used as a reference for computing $J_\fid$ in the following sections, where the wolf suppressor(s) will be enabled. 

\begin{figure}[h!]
	\centering
	\includegraphics[width=\textwidth]{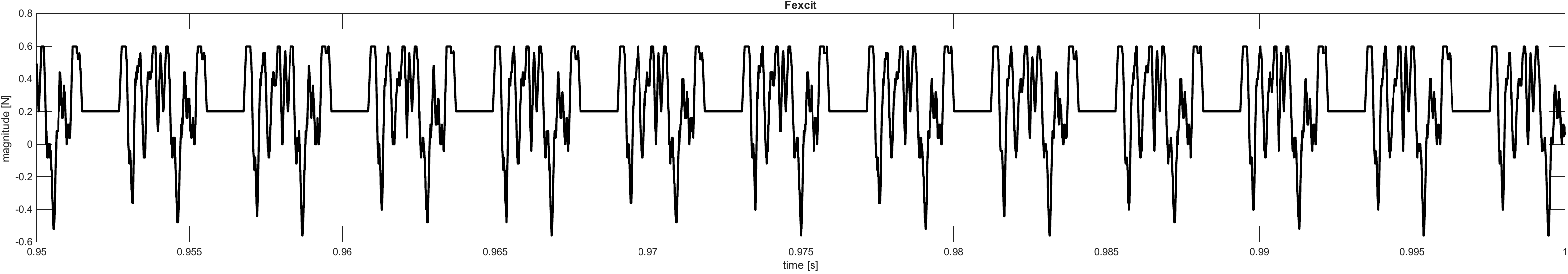} 
	\caption{$F_\bow$ as a function of time (zoom).}
	\label{fig:bow_0notch_f1_Fbow}
\end{figure}

Figure \ref{fig:bow_0notch_f1_Fbow} shows the bow force $F_\bow$ defined in \eqref{def:Fbow} when the note $f_1$ is being played. The force is slightly smoothed via a moving average for better readability.
We can observe the regular alternation between the stick and slip phases.

\begin{figure}[h!]
	\centering
	\begin{subfigure}{0.49\textwidth}
		\centering
		\includegraphics[width=\textwidth]{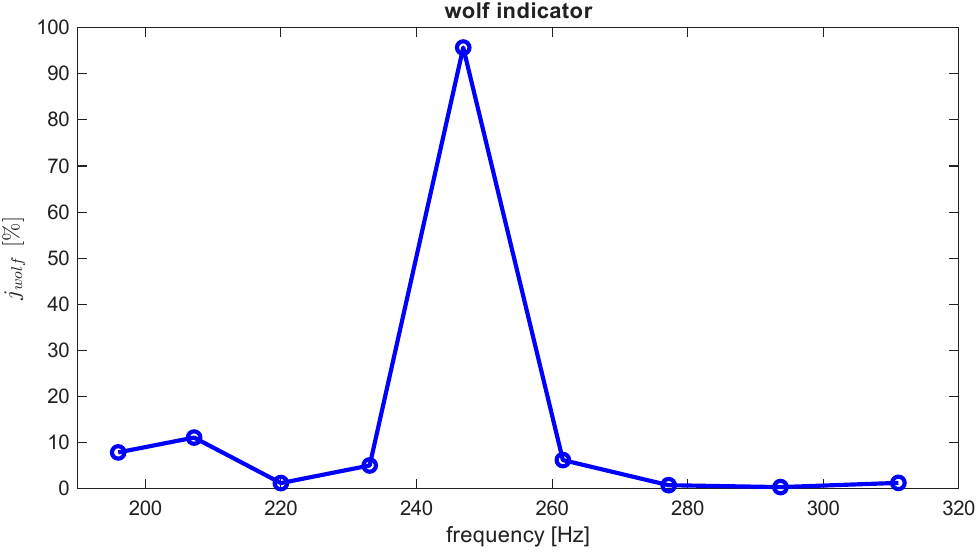} 
		\caption{Indicator $f_i\to j^i_\wolf$}
		\label{fig:bow_0notch_allf_Jwolf}
	\end{subfigure}
	\\
	\begin{subfigure}{0.49\textwidth}
		\centering
		\includegraphics[width=\textwidth]{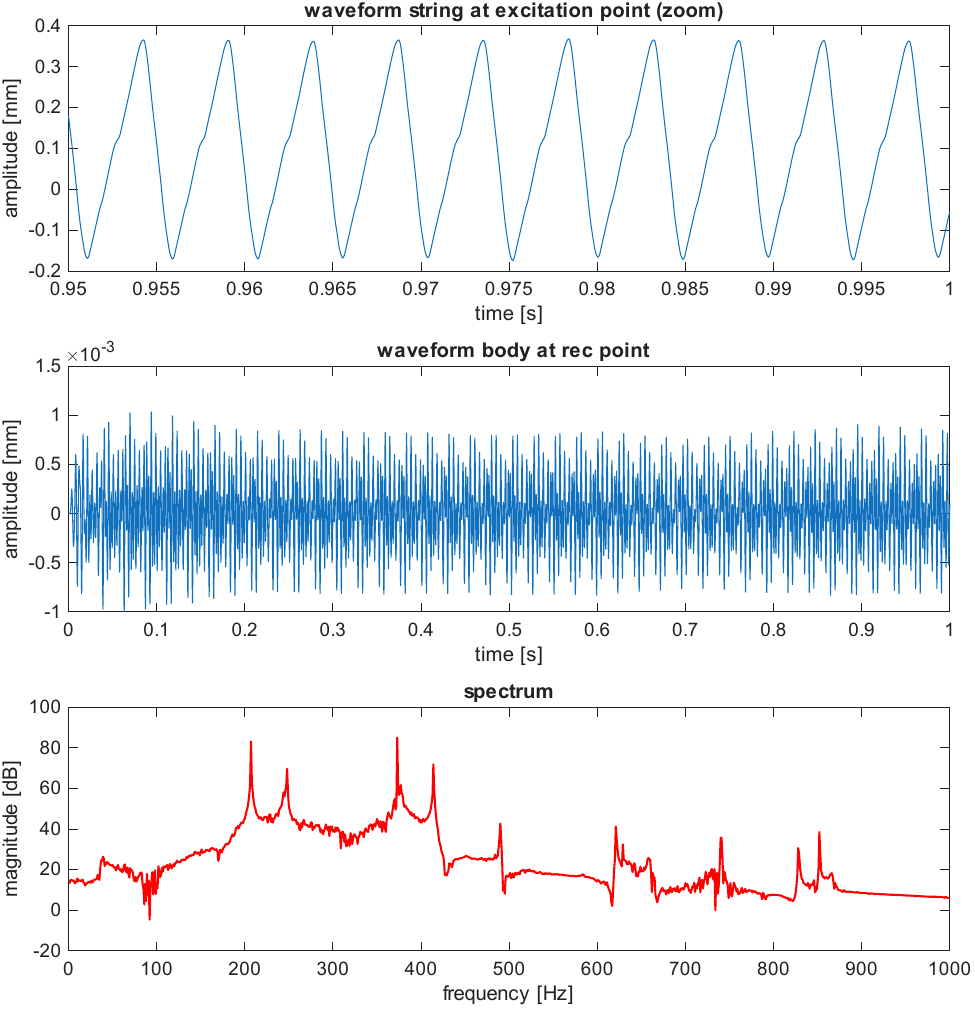}
		\caption{
			Waveform and spectrum for $f_2$ 	
                                       \audionotelocal{./supplementary_material/audio/9n_bow_0notch_f2_audio.wav}
			\audionoteonline{www.emilianocristiani.it/attach/paper_wolfnote/audio/9n_bow_0notch_f2_audio.wav}
		}
		\label{fig:bow_0notch_f2_WF_SP}
	\end{subfigure}
	\begin{subfigure}{0.49\textwidth}
		\centering
		\includegraphics[width=\textwidth]{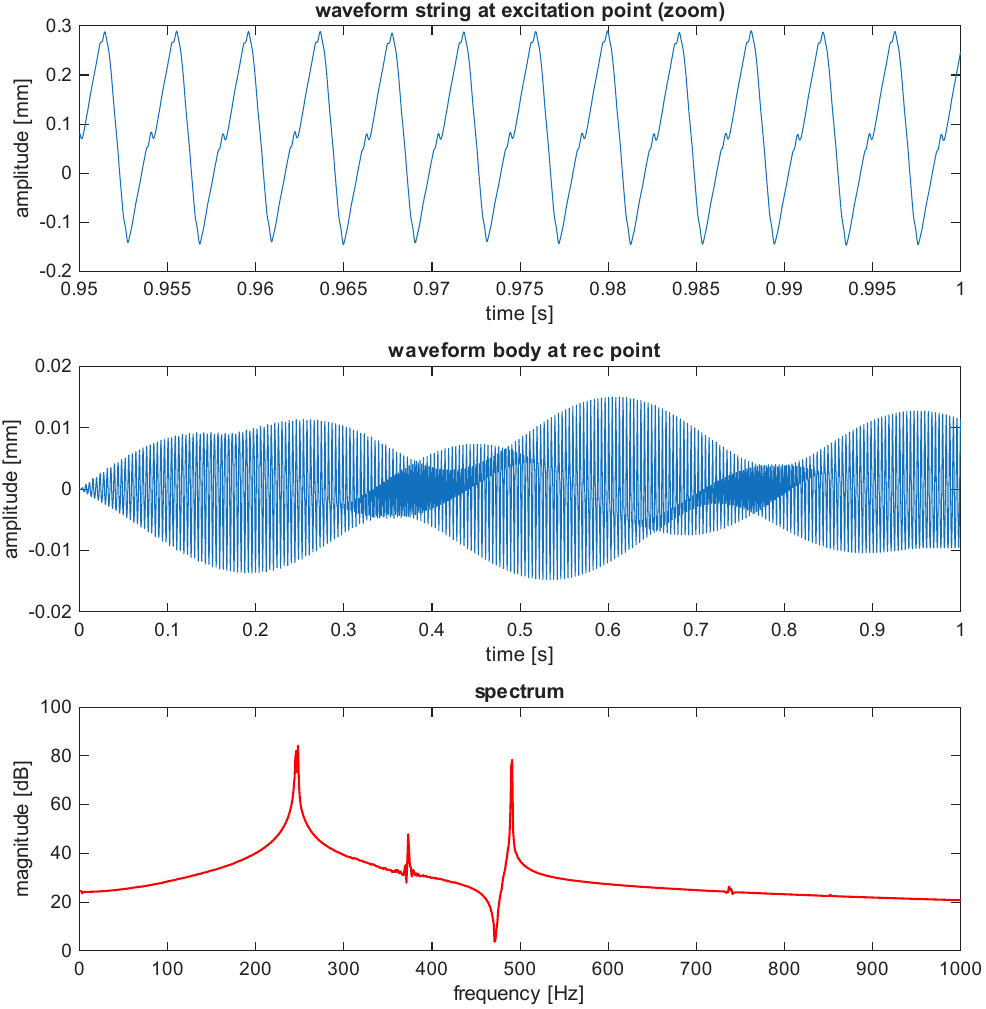}
		\caption{
			Waveform and spectrum for $f_5$ 
                                       \audionotelocal{./supplementary_material/audio/9n_bow_0notch_f5_audio.wav}
			\audionoteonline{www.emilianocristiani.it/attach/paper_wolfnote/audio/9n_bow_0notch_f5_audio.wav}
		}
		\label{fig:bow_0notch_f5_WF_SP}
	\end{subfigure}
	\caption{
		BOW-0S: Indicator $f_i\to j^i_\wolf$ and results for two frequencies from Table \ref{tab:frequencies}. 
		Waveforms of all frequencies are available 
		                               \href{./supplementary_material/large_figures/9n_bow_0notch_allf_WF.pdf}{here-local} and
		\href{http://www.emilianocristiani.it/attach/paper_wolfnote/large_figures/9n_bow_0notch_allf_WF.pdf}{here-online}.
		Sounds of all frequencies are available 
		                               \href{./supplementary_material/audio/9n_bow_0notch_allf_audio.wav}{here-local} and
		\href{http://www.emilianocristiani.it/attach/paper_wolfnote/audio/9n_bow_0notch_allf_audio.wav}{here-online}.
	}  
	\label{fig:bow_0notch}
\end{figure}

Figure \ref{fig:bow_0notch_allf_Jwolf} shows the indicator $j^i_\wolf$, defined in \eqref{def:jiwolf}, for any frequency $i=1,\ldots,N_f$. 
We see that again $j^5_\wolf>95\%$, meaning that the wolf note is captured with very high precision. 

Figures \ref{fig:bow_0notch_f2_WF_SP}-\ref{fig:bow_0notch_f5_WF_SP} show waveforms and spectra of $f_2$ and $f_5$, respectively. 
The wolf tone at $f_5$ appears again as in the pluck case, but it is less regular.

\clearpage
\subsubsection{Test BOW-1S: optimal placement of one wolf suppressor}
As we did in Section \ref{sec:testpluck-1S}, we enable one wolf suppressor and investigate the extent to which it is able to prevent the onset of the wolf note. 
We also study how the suppressor perturbs notes that were not originally affected by the wolf.

\begin{figure}[h!]
	\centering
	\begin{subfigure}{0.33\textwidth}
		\centering
		\includegraphics[width=\textwidth]{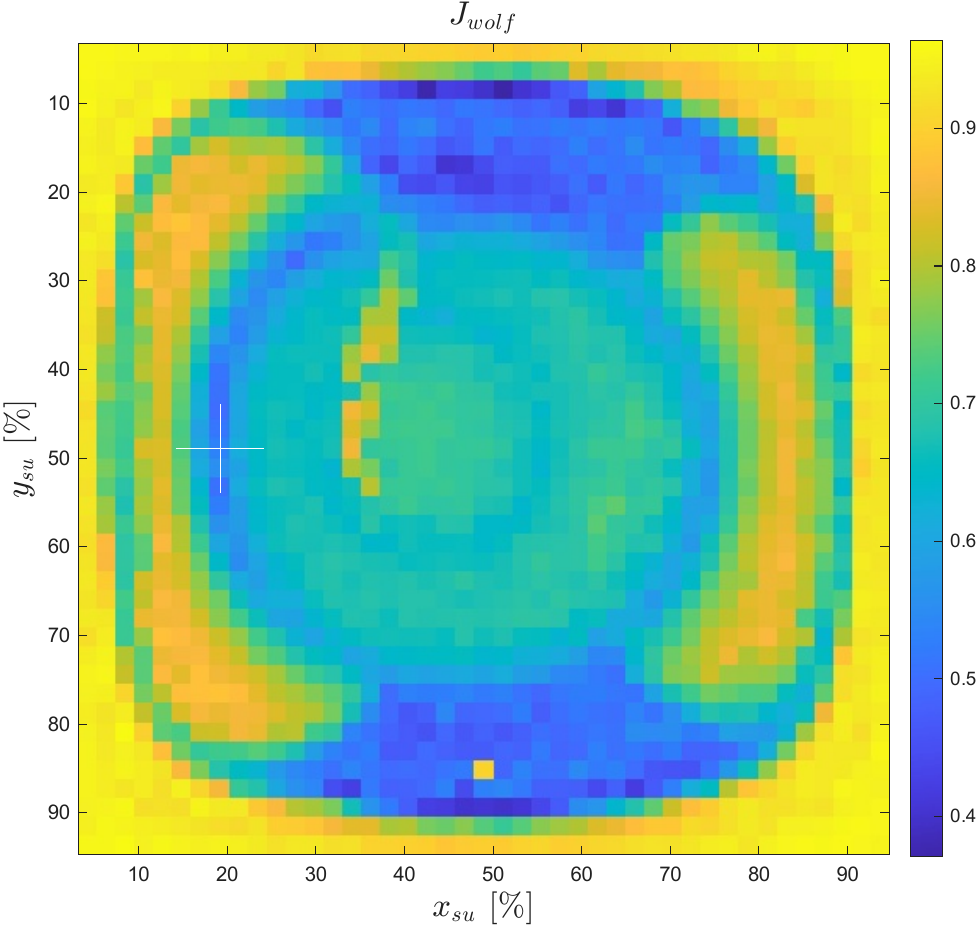} 
		\caption{$J_\wolf$}
		\label{fig:bow_1notch_allf_Jwolf_landscape}
	\end{subfigure}
	\begin{subfigure}{0.33\textwidth}
		\centering
		\includegraphics[width=\textwidth]{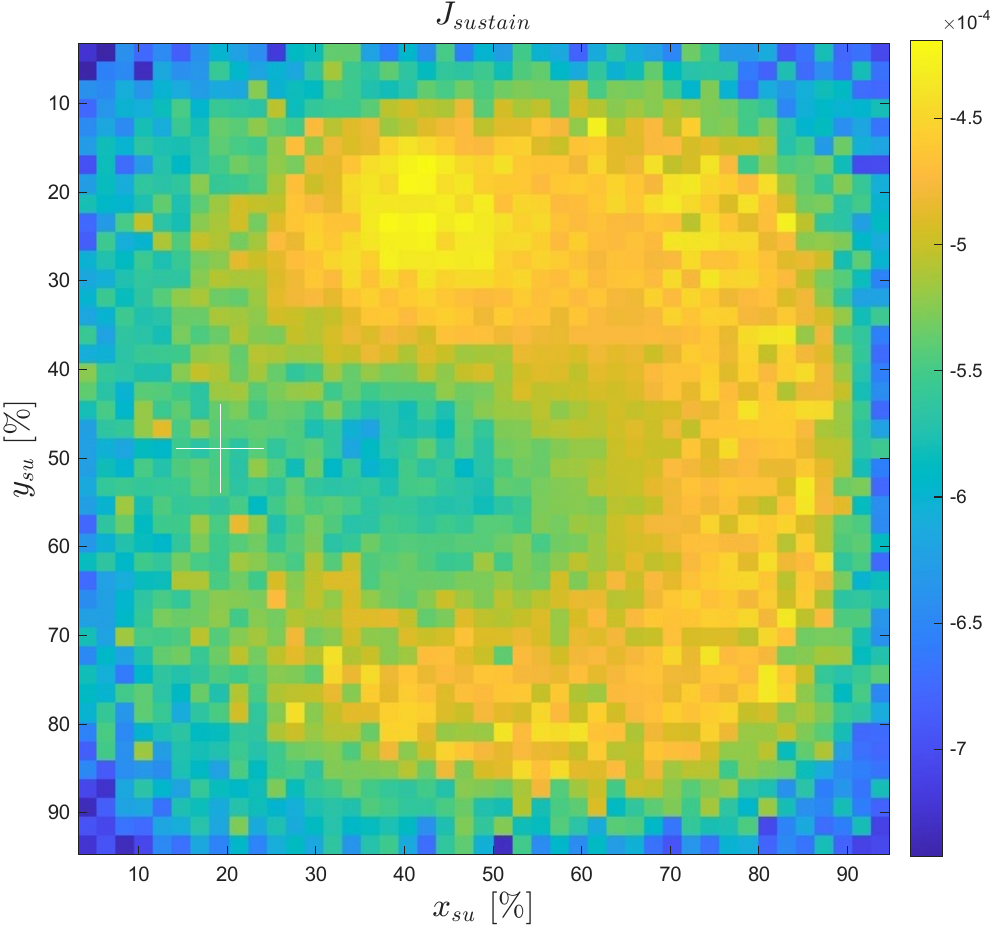} 
		\caption{$J_\clos$}
		\label{fig:bow_1notch_allf_Jclos_landscape}
	\end{subfigure}
    \begin{subfigure}{0.33\textwidth}
		\centering
		\includegraphics[width=\textwidth]{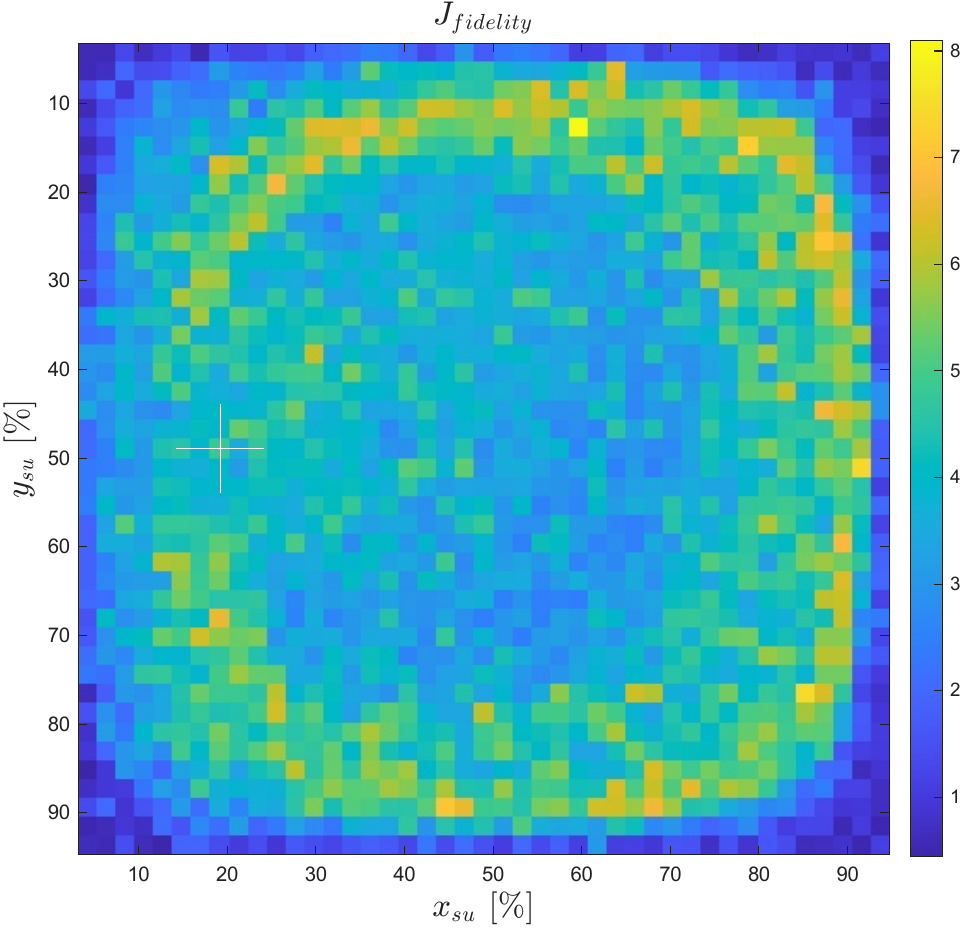} 
		\caption{$J_\fid$}
		\label{fig:bow_1notch_allf_Jfid_landscape}
	\end{subfigure}
	\caption{BOW-1S. Heat maps of $J_\wolf$, $J_\clos$, and $J_\fid$ as a function of the wolf suppressor position $(x_\su,y_\su)$ on the instrument body (expressed as $\%$ of the body size). The white cross indicates the point chosen as best placement of the suppressor.
	}
	\label{fig:bow_1notch_allf_allJ_landscape}
\end{figure}

Figure \ref{fig:bow_1notch_allf_allJ_landscape} shows the heat maps of the indicators $J_\wolf$, $J_\clos$, and $J_\fid$ as a function of the suppressor position $(x_\su,y_\su)$.  

We immediately observe that the images are different from those shown in Figure \ref{fig:pluck_1notch_allf_allJ_landscape}, this time being more chaotic although a structure is preserved. 

For the following tests, we heuristically choose the point $(19\%,49\%)$ (indicated with the white cross in the figure) as the best compromise. Indeed, a wolf suppressor located there is able to eliminate the wolf tone while still preserving a good sustain and, thirdly in importance, fidelity. 

\begin{figure}[h!]
	\centering
	\begin{subfigure}{0.49\textwidth}
		\centering
		\includegraphics[width=\textwidth]{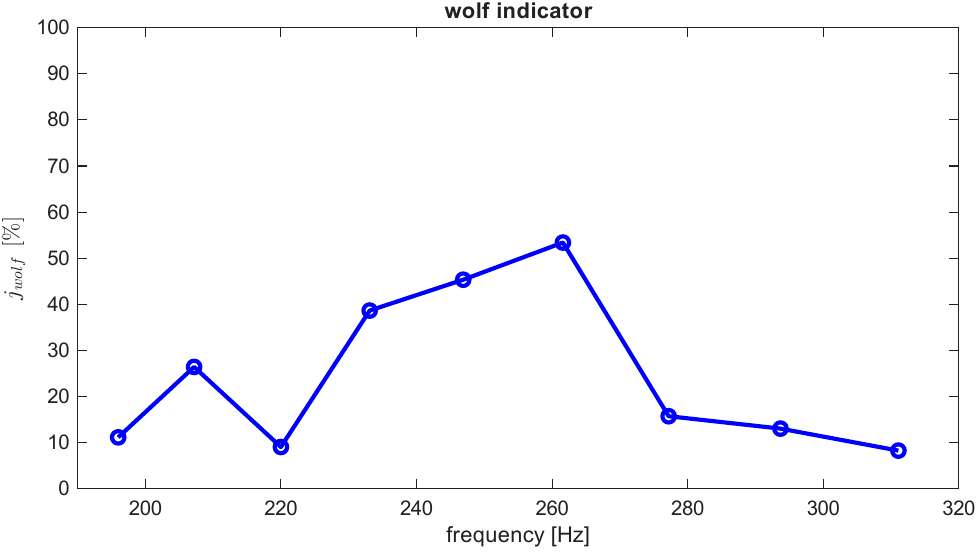} 
		\caption{Indicator $f_i\to j^i_\wolf$}
		\label{fig:bow_1notch0923_allf_Jwolf}
	\end{subfigure}
	\\
	\begin{subfigure}{0.49\textwidth}
		\centering
		\includegraphics[width=\textwidth]{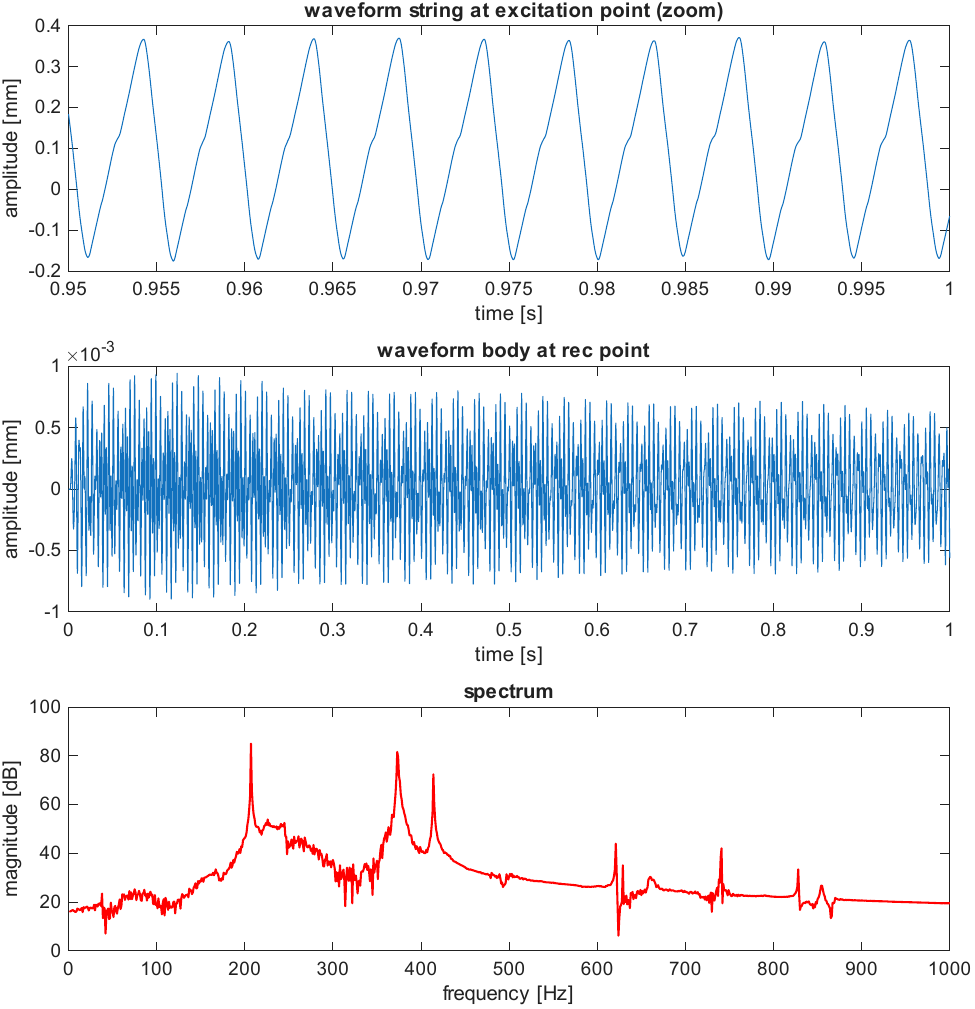}
		\caption{
			Waveform and spectrum for $f_2$ 	
			\audionotelocal{./supplementary_material/audio/9n_bow_1notch0923_f2_audio.wav}
			\audionoteonline{www.emilianocristiani.it/attach/paper_wolfnote/audio/9n_bow_1notch0923_f2_audio.wav}
		}
		\label{fig:bow_1notch0923_f2_WF_SP}
	\end{subfigure}
	\begin{subfigure}{0.49\textwidth}
		\centering
		\includegraphics[width=\textwidth]{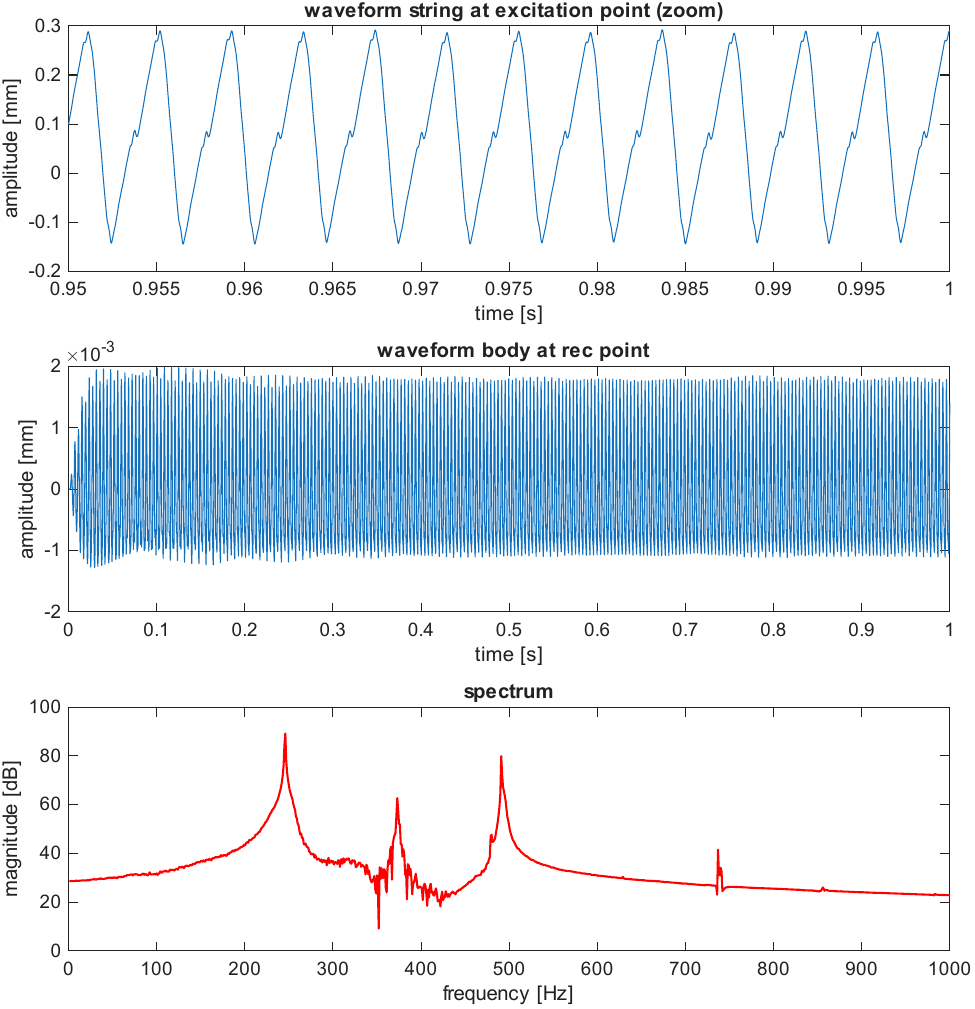}
		\caption{
			Waveform and spectrum for $f_5$ 
			\audionotelocal{./supplementary_material/audio/9n_bow_1notch0923_f5_audio.wav}
			\audionoteonline{www.emilianocristiani.it/attach/paper_wolfnote/audio/9n_bow_1notch0923_f5_audio.wav}
		}
		\label{fig:bow_1notch0923_f5_WF_SP}
	\end{subfigure}
	\caption{
		BOW-1S: Indicator $f_i\to j^i_\wolf$ and results for two frequencies from Table \ref{tab:frequencies}. 
		Waveforms of all frequencies are available 
		\href{./supplementary_material/large_figures/9n_bow_1notch0923_allf_WF.pdf}{here-local} and
		\href{http://www.emilianocristiani.it/attach/paper_wolfnote/large_figures/9n_bow_1notch0923_allf_WF.pdf}{here-online}.
		Sounds of all frequencies are available 
		\href{./supplementary_material/audio/9n_bow_1notch0923_allf_audio.wav}{here-local} and
		\href{http://www.emilianocristiani.it/attach/paper_wolfnote/audio/9n_bow_1notch0923_allf_audio.wav}{here-online}.
	}  
	\label{fig:bow_1notch}
\end{figure}

Figure \ref{fig:bow_1notch0923_allf_Jwolf} shows the indicator $j^i_\wolf$ for any frequency $i=1,\ldots,N_f$. We see that all frequencies are below 60\%, meaning that the wolf note is well suppressed. 

Figures \ref{fig:bow_1notch0923_f2_WF_SP}-\ref{fig:bow_1notch0923_f5_WF_SP} show waveforms and spectra of $f_2$ and $f_{5}$, respectively. 
Comparing the results with those from Section \ref{sec:testbow-0S}, we confirm that the wolf tone is eliminated but some notes are attenuated.
%

\clearpage
\subsubsection{Test BOW-2S: optimal placement of two wolf suppressors}\label{sec:testbow-2S}
In this section, we try to add a second wolf suppressor. 
To be fair, we choose two suppressors each weighing half of the one used in the single case.
By doing this, we aim at understanding whether two suppressors with the same total mass can be somehow better than one.
Due to the high computational cost, we do not run an exhaustive computation that visits all possible placements of both suppressors; instead, we move only one coordinate of each suppressor. 
More precisely, we fix $x_{\su_1}=42\%$ for the first suppressor and $y_{\su_2}=50\%$ for the second one, and let the other coordinates vary from 0 to 100\% as before.

\begin{figure}[h!]
	\centering
	\begin{subfigure}{0.33\textwidth}
		\centering
		\includegraphics[width=\textwidth]{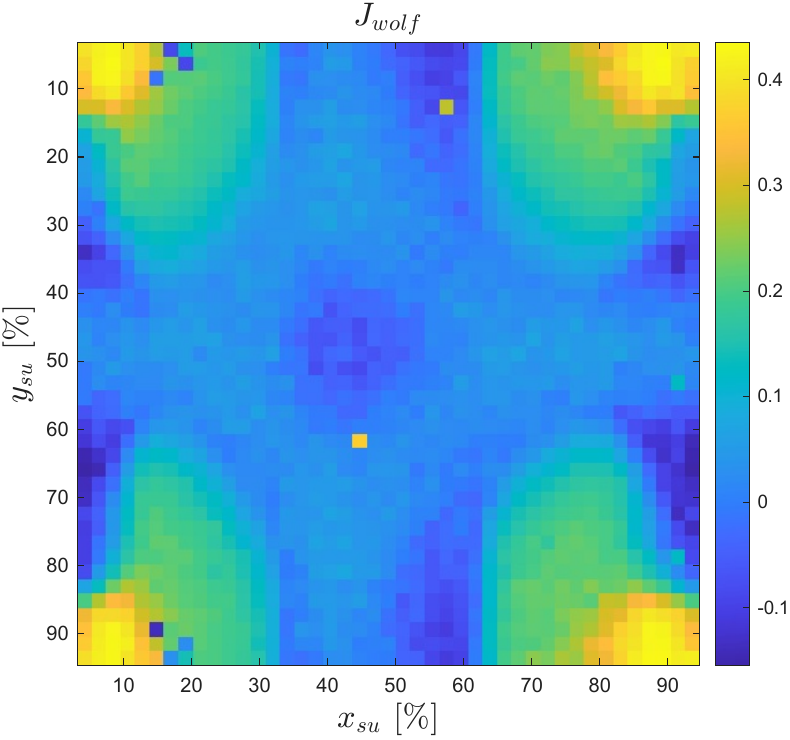} 
		\caption{$J_\wolf$}
		\label{fig:bow_2notch_cross_allf_Jwolf_landscape}
	\end{subfigure}
	\begin{subfigure}{0.33\textwidth}
		\centering
		\includegraphics[width=\textwidth]{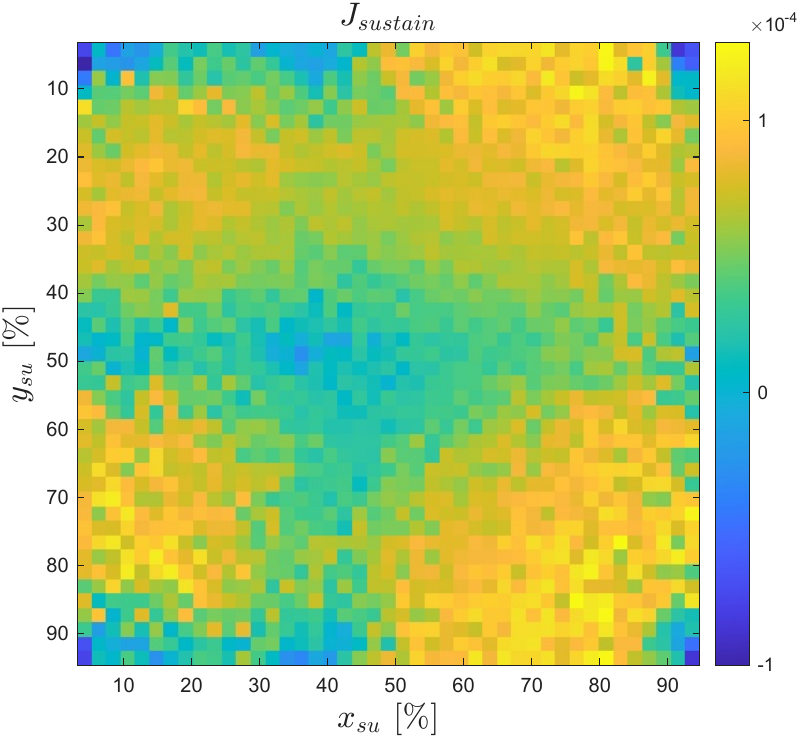} 
		\caption{$J_\clos$}
		\label{fig:bow_2notch_cross_allf_Jclos_landscape}
	\end{subfigure}
    \begin{subfigure}{0.33\textwidth}
		\centering
		\includegraphics[width=\textwidth]{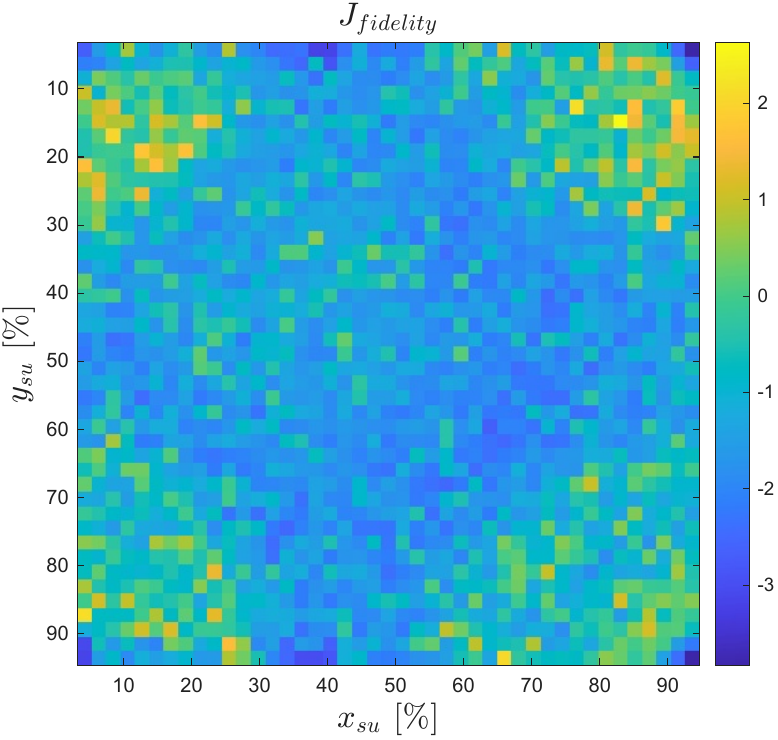} 
		\caption{$J_\fid$}
		\label{fig:bow_2notch_cross_allf_Jfid_landscape}
	\end{subfigure}
	\caption{BOW-2S. Heat maps of $J_\wolf$, $J_\clos$, and $J_\fid$ as a function of the $x$ coordinate of first suppressor and the $y$ coordinate of the second suppressor varying on the instrument body (expressed as $\%$ of the body size). 
	The colorbar's scale is defined in such a way that values, when negative, quantify the \emph{improvement} with respect to the best values which can be obtained with only one wolf suppressor.
	}
	\label{fig:bow_2notch_cross_allf_allJ_landscape}
\end{figure}

Figure \ref{fig:bow_2notch_cross_allf_allJ_landscape} shows the relative variation of the three indicators with respect to the best one-suppressor configuration.

We observe that we are able to find positions that allow us to improve the indicators with respect to the single wolf suppressor usage. Unfortunately, the gain is limited and is not actually perceptible.
%
%
%
%
%
%
\section{Conclusions and future work}
We have proposed and analyzed a coupled string-bridge-body-suppressors model to investigate the wolf note phenomenon and the effect of passive suppressors, introducing quantitative indicators to measure wolf tone, attenuation, and spectral fidelity. 
The proposed $J_\wolf$ indicator was also tested on a real cello recording, where it correctly identifies the wolf note (see Appendix).

Numerical results highlight that bowed excitation produces more irregular dynamics than plucking, making wolf suppression more challenging in realistic conditions. 
Moreover, while adding multiple suppressors can further reduce the wolf effect and the other indicators, the improvement is marginal and comes at the cost of a more complex optimal placement strategy.

Future work will focus on extending the model toward a fully three-dimensional body including air coupling, and developing multi-objective optimization strategies. 

\bibliographystyle{emibibliostyle}
\bibliography{bibliomusic}

\appendix
\section{Real-life wolf note}
\begin{figure}[h!]
	\centering
	\includegraphics[width=\textwidth]{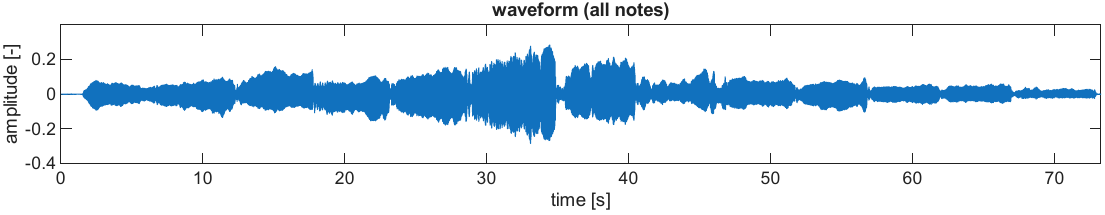}
	\caption{
		Waveform for chromatic scale C3-C4 (13 notes).
		\audionotelocal{./supplementary_material/audio/real_Do3Do4.mp4} 	
		\audionoteonline{www.emilianocristiani.it/attach/paper_wolfnote/audio/real_Do3Do4.mp4}
	}
	\label{fig:real_all_notes}
\end{figure}

\begin{figure}[h!]
	\begin{subfigure}{0.49\textwidth}
		\centering
		\includegraphics[width=\textwidth]{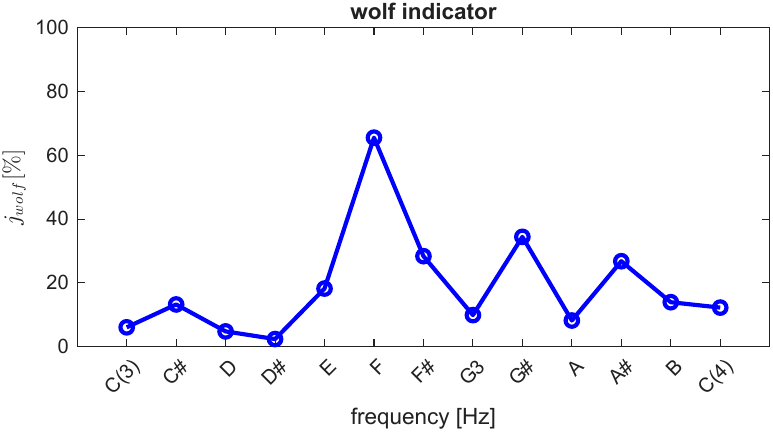}
		\caption{Indicator $f_i\to j^i_\wolf$}
		\label{fig:realwolf_jwolf}
	\end{subfigure}
	\begin{subfigure}{0.49\textwidth}
		\centering
		\includegraphics[width=0.96\textwidth]{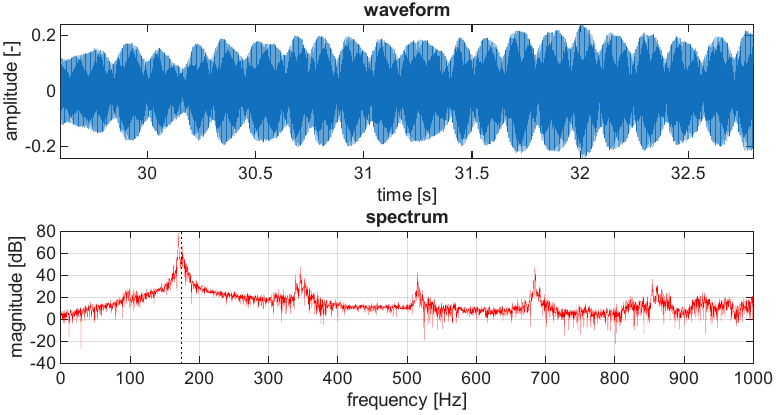}
		\caption{Waveform and spectrum for the wolf note F. The black dotted vertical line represents the ideal F frequency}
		\label{fig:realwolf_WF_SP}
	\end{subfigure}
	\caption{Indicator $f_i\to j^i_\wolf$ and analysis of the real-life wolf note.}
	\label{fig:realwolf}
\end{figure}

Here we briefly analyze a recording from a real cello playing a chromatic scale from C3 to C4. Figure \ref{fig:real_all_notes} shows the waveform for the entire chromatic scale. 
Figure \ref{fig:realwolf_jwolf} reports the corresponding values of the wolf indicator $j_\wolf$, defined in \eqref{def:jiwolf}, computed with the same parameters used throughout the paper.
One can observe that the note F, which exhibits the wolf, is well identified. 
Finally, Figure \ref{fig:realwolf_WF_SP} shows the waveform and the spectrum of the note F.

%
%
%
%
%
%

\end{document}